\documentclass[times,sort&compress,3p]{elsarticle}
\journal{Journal of Multivariate Analysis}

\usepackage[labelfont=bf]{caption}

\usepackage{amsmath,amsfonts,amssymb,amsthm,booktabs,color,epsfig,graphicx,hyperref,url,xspace,tikz,enumerate}

\usepackage{algorithmic}
\usepackage{subfig}
\usepackage[ruled]{algorithm2e}
\usepackage{stackrel}
\newcommand{\TPtwo}{$\mathrm{TP}_2$\xspace}
\newcommand{\PQD}{$\mathrm{PQD}$\xspace}
\newcommand{\SI}{$\mathrm{SI}$\xspace}
\newcommand{\SD}{$\mathrm{SD}$\xspace}
\newcommand{\LTI}{$\mathrm{LTI}$\xspace}

\newcommand{\LTD}{$\mathrm{LTD}$\xspace}

\newcommand{\NQD}{$\mathrm{NQD}$\xspace}
\newcommand{\Id}{\mathrm{Id}\xspace}
\newcommand{\LC}{Liebscher\xspace}
\newcommand{\FC}{{ comonotonic}\xspace}
\newcommand{\LF}{{comonotonic-based Liebscher}\xspace}
\newcommand{\CL}{\tilde{C}_{\mathrm{L}}}
\newcommand{\CF}{C_{\mathrm{F}}}
\newcommand{\CLF}{ \tilde{C}_{\mathrm{CL}}}
\newcommand{\SLF}{ \tilde{S}_{\mathrm{CL}}}
\newcommand{\ALF}{ \tilde{A}_{\mathrm{CL}}}
\newcommand{\K}{\mathcal{K}\xspace}
\newcommand{\obs}{\mathrm{obs}}

\newcommand{\algocomment}[1]{\hfill $\triangleright$ #1}

\usepackage{def_ju}

\newcommand{\CQFD}
{%
\mbox{}%
\nolinebreak%
\hfill%
\qed
\medbreak%
\par%
}

\newtheorem{Theo}{Theorem}
\newtheorem{Prop}{Proposition}
\newtheorem{Coro}{Corollary}
\newtheorem{Lem}{Lemma}
\newtheorem{Example}{Example}

\theoremstyle{remark}

\usepackage{bbold}
\usepackage{bm}

\begin{document}

\begin{frontmatter}

\title{Dependence properties and Bayesian inference for asymmetric multivariate copulas}

\author[A1]{Julyan Arbel}
\author[A1]{Marta Crispino}
\author[A1]{St\'ephane Girard\corref{mycorrespondingauthor}}

\address[A1]{Univ. Grenoble Alpes, Inria, CNRS, Grenoble INP, LJK, 38000 Grenoble, France}

\cortext[mycorrespondingauthor]{Corresponding author. Email address: \url{stephane.girard@inria.fr}}

\begin{abstract}
We study  a broad class of asymmetric copulas introduced by \citet{Liebscher} as a combination of multiple---usually symmetric---copulas. 
The main thrust of the paper is to provide new theoretical properties including exact tail dependence expressions and stability properties. 
A subclass of Liebscher copulas obtained by combining comonotonic copulas is studied in more details. 
We establish further dependence properties for copulas of this class and show that they are characterized by an arbitrary number of singular components. 
Furthermore, we introduce a novel iterative representation for general Liebscher copulas which \textit{de facto} insures uniform margins, thus relaxing a constraint of Liebscher's original construction. 
Besides, we show that this iterative construction proves useful for  inference by developing an Approximate Bayesian computation sampling scheme.  
This inferential procedure is demonstrated on simulated data and is compared to a likelihood-based approach in a setting where the latter is available.
\end{abstract}


\end{frontmatter}

\section{Introduction\label{sec:1}}

Let $\bm X= (X_1,\ldots ,X_d )$ be a continuous random vector with $d$-variate cumulative distribution function (cdf) $F$, and let $F_j$, $j\in\{1,\ldots,d\}$, be the marginal cdf of $X_j$. According to Sklar's theorem~\citep{sklar1959}, there exists a unique $d$-variate function $C: [0,1]^d\rightarrow[0,1]$ such that
$$F(\bm{x}) = C(F_1(x_1),\ldots,F_d(x_d)), \quad \bm{x}=(x_1,\ldots,x_d) \in\mathbb{R}^d.$$
The function $C$ is referred to as the copula associated with $F$. It is the $d$-dimensional cdf of the random vector $(F_1(X_1),\ldots,F_d(X_d))$ with uniform margins on $[0,1]$. 

A copula is said to be symmetric (or exchangeable) if for any $\bm{u}\in[0,1]^d$, and for any permutation $(\sigma_1,\ldots,\sigma_d)$ of the first $d$ integers $\{1,\ldots,d\}$, it holds that $C(u_1,\ldots,u_d)=C(u_{\sigma_1},\ldots,u_{\sigma_{d}}).$
The assumption of exchangeability may be unrealistic in many domains, including quantitative risk management \citep{di2016asymmetric}, reliability modeling \citep{wu2014construction}, and oceanography \citep{zhang2018modeling}. The urge for asymmetric copula models in order to better account for complex dependence structures has recently stimulated research in several directions, including \citet{rodriguez2004new, alfonsi2005new, durante2009construction, wu2014construction, durante2015copulas}. We focus here on a simple yet general method for building asymmetric copulas introduced by~\citet[Theorem 2.1,  Property (i)]{Liebscher}:

\begin{Theo}\citep[Liebscher, ][]{Liebscher}
\label{theoLieb}
Let $C_1,\ldots,C_K : [0,1]^d \to [0,1]$ be copulas, $g_{j}^{(k)} : [0,1] \to [0,1]$ be increasing functions  such that $g_{j}^{(k)}(0)=0$ and $g_{j}^{(k)}(1)=1$ \text{for all} $k \in\{1,\ldots,K\}$ and $j \in\{1,\ldots,d\}$. 
Then,
\begin{equation}\label{eq:lieb_cop}
\bm{u}\in[0,1]^d\mapsto \tilde{C}(\bm{u})= \prod_{k=1}^K C_k(g_{1}^{(k)}(u_1),\ldots,g_{d}^{(k)}(u_d))
\end{equation}
is also a copula under the constraint that 
\begin{equation}\label{eq:lieb_assumption}
 \prod_{k=1}^K g_{j}^{(k)}(u) = u\quad  \text{ for all } u \in [0,1], j \in\{1,\ldots,d\}.
\end{equation}
\end{Theo}
Theorem~\ref{theoLieb} provides a generic way to construct an asymmetric copula $\tilde{C}$, henceforth referred to as \textit{Liebscher copula}, starting from a sequence of symmetric copulas $C_1,\ldots,C_K$. 
This mechanism was first introduced by~\cite{Khoudraji} in the particular case where $K=2$ and with the functions $g_j^{(k)}$ assumed to be power functions, for each $j\in\{1,\ldots,d\}$ and $k\in\{1,\ldots,K\}$, that satisfy condition~\eqref{eq:lieb_assumption}.
The class of Liebscher copulas covers a broad range of dependencies and benefits from tractable bounds on dependence coefficients of the bivariate marginals~\citep{Liebscher,Liebscher2,mazo2015class}. 
However, there are two main reasons why the practical implementation of this approach is not straightforward: (i) it is not immediate to construct functions that satisfy condition~\eqref{eq:lieb_assumption};
and (ii) the product form complicates the density computation even numerically, which makes it difficult to perform likelihood  inference on the model parameters~\citep{mazo2015class}. 

The aim of this paper is to deepen the understanding of \LC's construction in order to  overcome  drawbacks (i) and (ii). Our contributions in this regard are three-fold. 
First, we provide theoretical properties of the asymmetric copulas in \eqref{eq:lieb_cop}, including exact expressions of tail dependence  indices, thus complementing the partial results of \citet{Liebscher,Liebscher2}. 
Second, we give an iterative representation of \eqref{eq:lieb_cop}  which  has  the  advantage  to  relax assumption \eqref{eq:lieb_assumption} by automatically satisfying it. 
Third, we develop an inferential procedure and a sampling scheme that rely on the newly developed iterative representation.

The Bayesian paradigm proves very useful for inference in our context as it overcomes the problematic computation of the maximum likelihood estimate, which requires the maximization of a very complicated likelihood function (see recent contributions \citep{valle2018bayesian,ning2018nonparametric}). 
General Bayesian sampling solutions in the form of Markov chain Monte Carlo are not particularly well-suited neither since they  require the evaluation of that complex likelihood.   
Instead, we resort to Approximate Bayesian computation (ABC), a technique dedicated to models with complicated, or intractable, likelihoods (see \citep{abc,robert2018abc,karabatsos2018approximate} for recent reviews). ABC requires the ability to sample from the model, which is straightforward with our iterative representation of Liebscher  copula. 
The adequacy of ABC for inference in copula models was leveraged by 
\citet{grazian2017approximate}, although in the different context of empirical likelihood estimation. A reversed approach to ours is followed by \citet{li2017extending}, who make use of copulas in order to adapt ABC to high-dimensional settings.

Since its introduction, the construction by \LC has received much attention in the copula literature (e.g., \citep{salvadori2010multivariate,
durante2010construction,
lauterbach2015some}). However, most studies have been limited to simple cases where the product in~\eqref{eq:lieb_cop} has only two terms.  
We hope that our paper will contribute to the further spreading of \LC's copulas, because it allows to exploit their full potential by: (i) better understanding their properties; (ii) providing a novel construction, which facilitates their use with an arbitrary number $K$ of terms in~\eqref{eq:lieb_cop}; and (iii) giving a strategy to make inference on them. 

On top of what has been presented above, an additional contribution of this paper is to derive specific results for the subclass of \LC's copula when two or more comonotonic copulas are combined, which we call \LF copula.
This subclass is characterized by an arbitrary number of singular components. To the best of our knowledge, this is the first paper to investigate this copula's properties and to provide an inference procedure.\\

The paper is organized as follows. 
Section~\ref{sec:properties} provides some theoretical results concerning the properties of asymmetric \LC copulas, also presenting the novel iterative construction.
In Section~\ref{sec:liebFre}, we introduce and analyze the \LF copula.
Section~\ref{sec:simu} is dedicated to the inference strategy. It  demonstrates our approach on simulated data   and provides a comparison with a likelihood-based approach for a class of Liebscher copulas where maximum likelihood estimation is feasible.
We conclude with a short discussion in Section~\ref{sec:concl}. Proofs are postponed to the Appendix.

\section{Properties of the copula}\label{sec:properties}

In this section, some new properties of the copula $\tilde{C}$ are established, complementing the ones in~\citet{Liebscher,Liebscher2}.
Sections~\ref{subsec:tailprop} and~\ref{subsec:dependence} are dedicated to (tail) dependence properties. For the sake of simplicity, we focus on the case $d=2$ of bivariate copulas. Some stability properties of Liebscher's construction are highlighted in Section~\ref{subsec:misc}. 
Finally, an alternative construction to Liebscher copula~(\ref{eq:lieb_cop}) is introduced in Section~\ref{sub:constr}.

\subsection{Tail dependence}\label{subsec:tailprop}

The \textit{lower and upper tail dependence functions}, denoted by $\Lambda_L(C;\cdot)$ and $\Lambda_U(C,\cdot)$
respectively, are defined for all $(x,y)\in \R_+^2$ by
\begin{align*}
\Lambda_L(C;x,y)  = \lim_{\varepsilon \to 0} \frac{C(\varepsilon x,\varepsilon y)}{\varepsilon},\quad \text{and} \quad
\Lambda_U(C;x,y)  = x+y+\lim_{\varepsilon \to 0} \frac{C(1-\varepsilon x,1-\varepsilon y)-1}{\varepsilon},
\end{align*}
where $C$ is a given bivariate copula, see for instance~\cite{joe2010tail}. Note that these limits exist under a bivariate regular variation assumption, see~\cite{resnick2013extreme}, Section~5.4.2 for details.
When they exist, these functions are homogeneous (\citep{joe2010tail}, Proposition~2.2), {\it i.e.}, for all $t\in(0,1]$ and $(x,y)\in \R_+^2$,
$\Lambda(C;t x, ty) = t \Lambda(C;x,y)$, where $\Lambda$ is equal to $\Lambda_L$ or $\Lambda_U$. The \textit{lower and upper tail dependence coefficients}, denoted by $\lambda_L(C)$ and $\lambda_U(C)$ respectively, 
are defined as the conditional probabilities that a random vector associated with a copula $C$ belongs to lower or upper tail orthants given that a univariate margin takes extreme values:
\begin{align*}
\lambda_L(C) = \lim_{u\to0} \frac{C(u,u)}{u},\quad
\lambda_U(C) = 2- \lim_{u\to1} \frac{C(u,u)-1}{u-1}.
\end{align*}
These coefficients can also be interpreted in terms of the tail dependence functions: 
$\lambda_L(C) = \Lambda_L(C; 1,1)$ and $\lambda_U(C) = \Lambda_U(C; 1,1)$.
Conversely, in view of the homogeneity property, the behavior of the tail dependence functions on the diagonal is determined by the tail dependence coefficients: 
$\Lambda_L(C; t,t)= \lambda_L(C) t$ and $\Lambda_U(C; t,t)=\lambda_U(C)t $ for all $t\in(0,1]$.
The tail dependence functions for Liebscher copula are provided by Proposition~\ref{prop:L} which, in view of the previous remarks, 
allows us to derive the tail dependence coefficients in Corollary~\ref{cor:lambda}. Some of these results rely on the
notion of (univariate) regular variation. Recall that a positive function $g$ is said to be regularly varying with index $\gamma$
if $g(xt)/g(x)\to t^\gamma$ as $x\to\infty$ for all $t>0$, see~\cite{bingham1989regular}.

\begin{Prop}\label{prop:L}
Let $(x,y)\in \R_+^2$ and consider $\tilde{C}$ the bivariate copula defined by~(\ref{eq:lieb_cop}) with $d=2$. 
\begin{description}
\item{\rm{(i)}} Lower tail, symmetric case.
Assume that $g^{(k)}_1=g^{(k)}_2$ is a regularly varying function with index $\gamma^{(k)}>0$ for all $k\in\{1,\dots,K\}$.
Then, 
$$
\Lambda_L(\tilde{C}; x,y) = 
     \prod_{k=1}^K \Lambda_L(C_k; x^{\gamma^{(k)}}, y^{\gamma^{(k)}})  
$$
and, necessarily, $\sum_{k=1}^K \gamma^{(k)}=1$.
\item {\rm{(ii)}} Lower tail, asymmetric case.
Suppose there exists $k_0\in\{1,\dots,K\}$ such that $g_1^{(k_0)}(\varepsilon)/g_2^{(k_0)}(\varepsilon)\to 0$ as
$\varepsilon\to 0$. Then,
$$
\Lambda_L(\tilde{C}; x,y) = 0.
$$
\item{\rm{(iii)}} Upper tail, general case. Assume that, for all $k\in\{1,\dots,K\}$, $g^{(k)}_1$ and $g^{(k)}_2$ are 
differentiable at 1, with derivative at 1 denoted by $\text{\emph{d}}_1^{(k)}$ and $\text{\emph{d}}_2^{(k)}$ respectively. Then 
$$
\Lambda_U(\tilde{C}; x,y) = \sum_{k=1}^K \Lambda_U(C_k; \text{\emph{d}}_1^{(k)}x,\text{\emph{d}}_2^{(k)} y) 
$$
and, necessarily,
$\sum_{k=1}^K \text{\emph{d}}_j^{(k)}= 1$, for $j\in\{1, 2\}$.
\item {\rm{(iv)}} Upper tail, particular case. If, in addition to \rm{(iii)}, $\text{\emph{d}}_1^{(k)}=\text{\emph{d}}_2^{(k)}\eqqcolon \text{\emph{d}}^{(k)}$ for all $k\in\{1,\dots,K\}$, then
$$
\Lambda_U(\tilde{C}; x,y)  =
 \sum_{k=1}^K \text{\emph{d}}^{(k)} \Lambda_U(C_k; x, y)
$$
and, necessarily, $\sum_{k=1}^K \text{\emph{d}}^{(k)}= 1$.
\end{description}
\end{Prop}

\noindent Let us note that the functions $g_j^{(k)}$ considered by~\citet{Liebscher} and indexed by (I-III) in his Section~2.1 all satisfy the assumptions of Proposition~\ref{prop:L}.
The following result complements~\cite[Proposition~2.3]{Liebscher} and~\cite[Proposition~0.1]{Liebscher2} which provide
bounds on the tail dependence coefficients. Here instead, explicit calculations are provided.

\begin{Coro}
\label{cor:lambda}
Let $\tilde{C}$ be the bivariate copula defined by~(\ref{eq:lieb_cop}) with $d=2$.
\begin{description}
\item{\rm{(i)}} Lower tail, symmetric case. Under the assumptions of Proposition~\ref{prop:L}\rm{(i)},
$
\lambda_L(\tilde{C}) =  \prod_{k=1}^K \lambda_L(C_k) .
$
\item {\rm{(ii)}} Lower tail, asymmetric case. Under the assumptions of Proposition~\ref{prop:L}\rm{(ii)},
$
\lambda_L(\tilde{C}) = 0.
$
\item{\rm{(iii)}} Upper tail, general case. Under the assumptions of Proposition~\ref{prop:L}\rm{(iii)},
$
\lambda_U(\tilde{C}) = \sum_{k=1}^K \Lambda_U(C_k; \text{\emph{d}}_1^{(k)},\text{\emph{d}}_2^{(k)}) 
$
and, necessarily, 
$\sum_{k=1}^K \text{\emph{d}}_j^{(k)}= 1$, for $j\in\{1, 2\}$.
\item {\rm{(iv)}} Upper tail, particular case. Under the assumptions of Proposition~\ref{prop:L}\rm{(iv)},
$
\lambda_U(\tilde{C})  =
 \sum_{k=1}^K \text{\emph{d}}^{(k)} \lambda_U(C_k)
$
and, necessarily, 
$\sum_{k=1}^K \text{\emph{d}}^{(k)}= 1$.
\end{description}
\end{Coro}

\noindent It appears that the lower and upper tail dependence coefficients have very different behaviors.
In the case \rm{(i)} of a symmetric copulas, the lower tail dependence coefficient $\lambda_L$ is the product of the lower tail dependence coefficients associated with the components. Besides, $\lambda_L=0$ as soon as a
component $k_0$ has functions $g_1^{(k_0)}$ and $g_2^{(k_0)}$ with different behaviors at the origin (case \rm{(ii)}).
At the opposite, the upper tail dependence coefficient does not vanish even though a
component $k_0$ has functions $g_1^{(k_0)}$ and $g_2^{(k_0)}$ with different behaviors at 1 (case \rm{(iii)}).
In the particular situation where all components $k\in\{1,\ldots,K\}$ have functions $g_1^{(k)}$ and $g_2^{(k)}$ with the same behavior at 1 (case \rm{(iv)}), $\lambda_U$ is 
 a convex combination of the upper tail dependence coefficients associated with the components.

\subsection{Dependence}\label{subsec:dependence}

Let $(X,Y)$ be a pair of random variables with continuous margins and associated copula $C$. 
\begin{itemize}
\item $X$ and $Y$ are said to be \textit{totally positive of order 2}, \TPtwo (see~\citep{joe1997multivariate}), if for all $x_1<y_1,\, x_2<y_2$,
\begin{equation*}
{\Pr}(X\leq x_1,Y\leq x_2){\Pr}(X\leq y_1,Y\leq y_2) \geq {\Pr}(X\leq x_1,Y\leq y_2){\Pr}(X\leq y_1,Y\leq x_2).
\end{equation*}
Since this can be equivalently written in terms of $C$, 
we will write in short that $C$ is \TPtwo. 
\item $X$ and $Y$ are said to be \textit{positively quadrant dependent} (PQD) if 
$$
{\Pr}(X\leq x, Y\leq y)\geq {\Pr}(X\leq x) {\Pr}( Y\leq y) \mbox{ for all }(x,y).
$$
Since this property can be characterized by the copula property $C\geq \Pi$ where $\Pi$ denotes the independence copula, see for instance~\cite[Paragraph~5.2.1]{Nelsen}, 
we shall write for short that $C$ is PQD. The \textit{negatively quadrant dependence} (NQD) property is similarly defined by  $C\leq \Pi$.
\item $X$ and $Y$ are said to be \textit{left-tail decreasing} (LTD) if 
\begin{equation}
\begin{aligned}
\label{eq:LTI-def}
&{\Pr}(X\leq x| Y\leq y) \text{  is a decreasing function of } y \text{ for all } x,\\
&{\Pr}(Y\leq y| X\leq x) \text{ is  a decreasing function of } x \text{ for all } y.
\end{aligned}
\end{equation}
From~\citet[Theorem~5.2.5]{Nelsen}, this property can be characterized by the copula properties
\begin{equation}
\begin{aligned}
\label{eq:LTI-charact}
C(u,v)/u \mbox{ is decreasing in } u \mbox{ for all } v\in[0,1], \\
C(u,v)/v \mbox{ is decreasing in } v \mbox{ for all } u\in[0,1], \end{aligned}
\end{equation}
and we shall thus write that $C$ is LTD. The \textit{left-tail increasing} property (LTI) 
is similarly defined by reversing the directions of the monotonicity in~\eqref{eq:LTI-def} and~\eqref{eq:LTI-charact}. 
\item $X$ and $Y$ are said to be \textit{stochastically increasing} (SI) if 
\begin{equation}
\begin{aligned}
\label{eq:SI-def}
&{\Pr}(X> x| Y= y) \text{  is an increasing function of } y \text{ for all } x,\\
&{\Pr}(Y> y| X= x) \text{ is an increasing function of } x \text{ for all } y.
\end{aligned}
\end{equation}
From~\citet[Corollary~5.2.11]{Nelsen}, this property can be characterized by the copula properties
\begin{equation}
\begin{aligned}
\label{eq:SI-charact}
C(u,v) \mbox{ is a concave function of } u \mbox{ for all } v\in[0,1], \\
C(u,v) \mbox{ is a concave function of } v \mbox{ for all } u\in[0,1], 
\end{aligned}
\end{equation}
and we shall thus write that $C$ is SI. The \textit{stochastically decreasing} (SD) property is similarly defined by replacing increasing by decreasing in~\eqref{eq:SI-def} and concave by convex in~\eqref{eq:SI-charact}. 
\end{itemize}

In the next proposition, we show that under mild conditions, the above dependence properties are preserved under Liebscher's construction, thus complementing \LTD and \TPtwo properties established in~\citet[Proposition~2.2]{Liebscher}.

\begin{Prop}\label{prop:long_list_of_dependence}
If copulas  $C_1,\ldots,C_K$  all satisfy any of the properties defined above, \TPtwo, \PQD, \NQD, \LTD,  \LTI, \SI or \SD, then the same is satisfied for the Liebscher copula $\tilde{C}$ defined in~\eqref{eq:lieb_cop}---for \SI (respectively \SD), the $g_j^{(k)}$ functions in Theorem~\ref{theoLieb} are  additionally required to be concave functions (respectively convex functions) and the copulas $C_k$ to be twice differentiable, $k\in\{1,\ldots,K\}$, $j\in\{1,\ldots,d\}$.
\end{Prop}

\subsection{Stability properties}\label{subsec:misc}

Let us focus on the situation where the functions $g_j^{(k)}$ of Theorem~\ref{theoLieb} are power functions:
for all $j\in\{1,\dots,d\}$,  $k\in\{1,\dots,K\}$ and $t\in[0,1]$,
let 
\begin{equation}
\label{power}
    g_j^{(k)}(t)=t^{p_j^{(k)}}, \quad p_j^{(k)}\in(0,1), \sum_{\ell=1}^K p_j^{(\ell)}=1.
\end{equation}
Recall that a copula $C_{\#}$ is said to be max-stable if for all 
integer $n\ge 1$ and $(u_1,\dots,u_d)\in [0,1]^d$:
$$
C_{\#}^n(u_1^{1/n},\dots,u_d^{1/n})=C_{\#}(u_1,\dots,u_d).
$$
From~\cite[Proposition~3]{frees1998understanding}, it is clear that associating max-stable copulas $C_k$ 
with power functions~(\ref{power}) in Liebscher construction~(\ref{eq:lieb_cop}) still yields a max-stable copula. The goal
of this paragraph is to investigate to what extent this result can be generalized.
Our first result establishes the stability of a family of Liebscher copulas built from homogeneous functions.
More specifically, each copula $C_k(\cdot)$ in~(\ref{eq:lieb_cop}) is rewritten as 
$C(\cdot  \;| \;\theta_{k})$ where
\begin{equation}
    \label{eq-new}
     C(\cdot \;| \;\theta_{k}) \coloneqq  \prod_{i=1}^m  \varphi_i^{\theta_{ik}}(\cdot),
\end{equation}
with $\theta_k=(\theta_{1k},\dots,\theta_{mk})^T$ and $\varphi_i:[0,1]^d\to[0,1]$, $i\in\{1,\dots,m\}$.
\begin{Prop}
 \label{propstable1}
For all $j\in\{1,\dots,d\}$ and $t\in[0,1]$,
let $g_j^{(k)}$ be given by~(\ref{power}) where $p_j^{(k)}=p^{(k)}$ for all $k\in\{1,\dots,K\}$.
 Let $m>0$ and for all $i\in\{1,\dots,m\}$ introduce
$\varphi_i:[0,1]^d\to[0,1]$ such that $\ln \mathop \circ \varphi_i \mathop \circ \exp$ is homogeneous of degree $\lambda_i$.
For all $k\in\{1,\ldots,K\}$, assume that $C(\cdot \;| \;\theta_{k})$ in~(\ref{eq-new})
is a copula for some $\theta_{k}\in {\mathbb R}^m$.
Then, copula~(\ref{eq:lieb_cop}) is given for all $K\geq 1$ by
$$
\tilde{C}^{(K)}(\cdot) = C(\cdot \;|\; \tilde \theta_{K}),
$$
with $\tilde\theta_K=(\tilde \theta_{1K},\dots,\tilde \theta_{mK})^T$, and for all $i\in\{1,\dots,m\}$,
$$
\tilde \theta_{iK}=\sum_{k=1}^K \theta_{ik} (p^{(k)})^{\lambda_i}.
$$
\end{Prop}

\begin{Example}[Gumbel-Barnett copula $C_k$]
\label{exGumbel}
Let $C_k$ be the Gumbel-Barnett copula~\citep[Table~4.1]{Nelsen}. It can 
be written as
$$
C_k({\bm u})=C({\bm u}\;|\;\theta_k)=\prod_{j=1}^d u_j \exp\left(\theta_k\prod_{j=1}^d\ln(1/u_j)\right) = \varphi_1^{\theta_{1k}}({\bm u})\varphi_2^{\theta_{2k}}({\bm u})
$$
with $\theta_{1k}=1$, $\theta_{2k}=\theta_k\geq 1$,
$$
\varphi_1({\bm u})= \prod_{j=1}^d u_j\; \mbox{ and } \;\varphi_2({\bm u})=\exp\left(\prod_{j=1}^d\ln(1/u_j)\right).
$$
It thus fulfills the assumptions of Proposition~\ref{propstable1} with $m=2$, $\lambda_1=1$ and $\lambda_2=d$.
\end{Example}

\begin{Example}[Extreme-value copula $C_k$]\label{exMaxStable}
Extreme-value copulas exactly correspond to max-stable copulas and are characterized by their
tail-dependence function $L$ as:
$$
C_{\#}(u_1,\dots,u_d)=\exp( - L(\ln(1/u_1),\dots,\ln(1/u_d))),
$$
where $L:{\mathbb R}_+^d\to {\mathbb R}_+$ is homogeneous of degree 1, see for instance~\citep{gudendorf_extreme-value_2010}.
It is thus clear that every max-stable copula $C_k$ fulfills the assumptions of Proposition~\ref{propstable1} with $m=1$, $\theta_{1k}=1$, $\lambda_1=1$ and $\ln \circ  \varphi_1\circ\exp({\bm t})=-L(-{\bm t})$, $t\in {\mathbb R}_-^d$. Moreover,
$$\tilde\theta_{1K} = \sum_{k=1}^K p^{(k)} = 1$$
and thus $\tilde{C}^{(K)}=C$ for all $K\geq 1$.
\end{Example}

\noindent  It appears that max-stable copulas can be considered as fixed-points of Liebscher's construction~(\ref{eq:lieb_cop}).
The next result shows that, under mild assumptions, they are the only copulas verifying this property.

\begin{Prop}
 \label{propstable2}
For all $j\in\{1,\dots,d\}$, let $g_j^{(k)}$ be given by~(\ref{power}) where $p_j^{(k)}=p^{(k)}$ for all $k\in\{1,\dots,K\}$.
Assume $C_k=C$ for all $k\in\{1,\dots,K\}$ and let $\tilde{C}^{(K)}$ be the copula defined by~(\ref{eq:lieb_cop}).
Then, $\tilde{C}^{(K)}=C$ for all $K\geq 1$ and for all sequences  $p^{(1)}, \ldots, p^{(K)} \in (0,1)$  such that $\sum_{k=1}^K p^{(k)}=1$ if and only if $C$ 
is max-stable.
\end{Prop}

\noindent To complete the links with max-stable copulas, let us consider the situation where $C_k=C$ and $p_j^{(k)}=1/K$
for all $j\in\{1,\dots,d\}$ and $k\in\{1,\dots,K\}$. Liebscher's construction thus yields
$$
\tilde{C}^{(K)}({\bm u})\coloneqq \tilde{C}({\bm u})= C^K(u_1^{1/K},\dots,u_d^{1/K}),
$$
which is the normalized cdf associated with the maximum of $K$ independent uniform random vectors distributed according to the cdf $C$.
Therefore, as $K\to\infty$, $\tilde{C}^{(K)}$ converges to a max-stable copula under standard extreme-value assumptions on $C$.

\subsection{An iterative construction}\label{sub:constr}

Let  ${\cal F}$ be the class of increasing functions $f:[0,1]\to [0,1]$ such that $f(0)=0$, $f(1)=1$ and $\Id/f$ is increasing, where $\Id$ denotes the identity function.
For all $k\geq 1$ let $C_k$ be a $d$-variate copula and $f^{(k)}_j \in {\cal F}$ for all $j\in\{1,\dots,d\}$, with the assumption $f^{(1)}_j(t)=1$ for all $t\in[0,1]$. 
We propose the following iterative construction of copulas. For all $\bm{u}\in[0,1]^d$, consider the sequence defined  by
\begin{align}
 \label{def1}
 \tilde{C}^{(1)}(\bm{u})&= C_1(\bm{u}), \\
 \label{def2}
 \tilde{C}^{(k)}(\bm{u})&= C_k\left(\frac{u_1}{f^{(k)}_1(u_1)},\dots,\frac{u_d}{f^{(k)}_d(u_d)} \right)
 \tilde{C}^{(k-1)}\left(f^{(k)}_1(u_1),\dots,f^{(k)}_d(u_d)\right), \; k\geq 2.
\end{align}
Lemma~\ref{theocop} in \ref{appendix:aux} shows that $\tilde{C}^{(k)}$ is a $d$-variate copula, for all $k\geq 1$.

Let $K\geq 1$. For all functions $f^{(1)},\dots,f^{(K)}$ : $[0,1]\to[0,1]$ and $i,j\in\{1,\dots,K\}, $ let us introduce the notation 
\begin{equation}
    \label{notation-rond}
\bigodot_{k=i}^j f^{(k)} \coloneqq  f^{(i)} \mathop \circ \dots \mathop \circ f^{(j)}
\mbox{ if } i\leq j \mbox{ and }
\bigodot_{k=i}^j f^{(k)}  \coloneqq \Id \mbox{ otherwise.}
\end{equation}
The next result shows that there is a one-to-one correspondence between copulas built by the iterative procedure~(\ref{def1}),~(\ref{def2}) 
and Liebscher copulas, reported in Theorem~\ref{theoLieb}.
\begin{Prop}
 \label{prop-rewrite}
 The copula $\tilde{C}^{(K)}$, $K\geq 1$ defined iteratively by~(\ref{def1}),~(\ref{def2}) is a Liebscher copula. It can be rewritten as
 \begin{equation}
 \label{def3}
  \tilde{C}^{(K)}(\bm{u}) = \prod_{k=1}^K C_k\left(g_1^{(K-k+1,K)}(u_1),\dots,g_d^{(K-k+1,K)}(u_d) \right)
 \end{equation}
 for all $\bm{u}\in[0,1]^d$
where, for all $j\in\{1,\dots,d\}$ and $K\geq 1$,
\begin{align}
\label{gj1}
g_j^{(1,K)}&=\Id/f_j^{(K)}, \\
\label{gj2}
g_j^{(k,K)}&= \bigodot_{i=K-k+2}^K f_j^{(i)} \left/ \bigodot_{i=K-k+1}^K f_j^{(i)} \right.,\, 
k\in\{2,\dots,K\}.
\end{align}
Conversely, each Liebscher copula (defined in Theorem~\ref{theoLieb}) can be built iteratively from~(\ref{def1}),~(\ref{def2}). 
\end{Prop}
\noindent Let us note that the iterative construction~(\ref{def1}),~(\ref{def2}) 
thus provides a way to build functions~(\ref{gj1}), (\ref{gj2}) that automatically fulfill Liebscher's  constraints~\eqref{eq:lieb_assumption} of Theorem~\ref{theoLieb}.
As a consequence, the construction~(\ref{def1}),~(\ref{def2}) also gives an iterative way to sample from a \LC copula~\eqref{eq:lieb_cop}, described in detailed in  Algorithm~\ref{algo:Iterative_copula_sampling} (see the proof of Lemma~\ref{theocop} in \ref{appendix:aux} for a theoretical justification).

\vspace{0.5cm}
\begin{center}
\begin{minipage}{10.8cm}
\begin{algorithm}[H]
\begin{footnotesize}
\DontPrintSemicolon
{\textbf{Input} $\left[f_j^{(k)}\right]_{k,j},(C_k)_{k}$}
\algocomment{functions in $\mathcal{F}$ appearing in~\eqref{def2}, and copulas}\;
{$(X_1^{(1)},\dots,X_d^{(1)})\sim {C}_1$}\;
\For{$k=2,\dots,K$}{
$(Y_1,\dots,Y_d)\sim {C}_k$ independently of $(X_1^{(k-1)},\dots,X_d^{(k-1)})$\;
\For{$j\in\{1,\dots,d\}$}{
$X_j^{(k)}=\max\left(\left(f^{(k)}_j\right)^{-1}(X_j^{(k-1)}),\left(\Id/f_j^{(k)}\right)^{-1}(Y_j)\right)$}}
{\textbf{Output} $\bm{X}=(X_1^{(K)},\ldots,X_d^{(K)})\sim \tilde{C}$}
\caption{Iterative sampling scheme for \LC copula~\eqref{eq:lieb_cop}}\label{algo:Iterative_copula_sampling}
\end{footnotesize}
\end{algorithm}
\end{minipage}
\end{center}
\vspace{0.5cm}

\begin{Example}[Power functions $f_j^{(k)}$]
\label{expower}
Let functions $f_j^{(k)}$ be power functions in the form of %
\begin{align}\label{eq:def_powers}
    f_j^{(k)}(t)=t^{1-a_j^{(k)}}, \quad a_j^{(1)}=1, \,\,    a_j^{(k)}\in(0,1), \,\, \text{ for all } k\geq 2,
\end{align}
for all $j\in\{1,\dots,d\}$ and $t\in[0,1]$.
From Proposition~\ref{prop-rewrite}, $g_j^{(k)}(t)=g_j^{(k,K)}(t)=t^{p_j^{(k,K)}}$ with 
\begin{equation}
\left\{
\begin{aligned}\label{eq:p-from-a}
p_j^{(1,K)} &= a_j^{(K)},\\
p_j^{(k,K)} &= a_j^{(K-k+1)} \prod_{i=K-k+2}^K (1-a_j^{(i)}), \; \mbox{ if } 2\leq k\leq K,
\end{aligned}
\right.
\end{equation}
for all $K\geq 1$ and  $j\in\{1,\dots,d\}$. Note that, by construction,
$$
\sum_{k=1}^K p_j^{(k,K)} =1,
$$
for all $j\in\{1,\dots,d\}$ and thus $\left(p_j^{(1,K)},\dots,p_j^{(K,K)}\right)$ can be interpreted as a discrete
probability distribution on $\{1,\dots,K\}$. Besides, let $\tilde{a}_j^{(k,K)}\coloneqq a_j^{(K+1-k)}$
for all $k\in\{1,\dots,K\}$. Equations~(\ref{eq:p-from-a}) can be rewritten as
$$
\left\{
\begin{aligned}
p_j^{(1,K)} &= \tilde{a}_j^{(1,K)},\\
p_j^{(k,K)} &= \tilde{a}_j^{(k,K)} \prod_{i=1}^{k-1} (1-\tilde{a}_j^{(i,K)}), \; \mbox{ if } 2\leq k\leq K,
\end{aligned}
\right.
$$
which corresponds to  the so-called stick-breaking construction  \citep{sethuraman1994constructive}. 
\end{Example}


\section[The \LF copula]{The \LF copula}\label{sec:liebFre}

We analyze here in more details the Liebscher copula obtained by combining $K\geq 2$ \FC (also called Fr\'echet) copulas defined by ${C}(u,v) = \min(u,v)$. We here focus on the bivariate case $d=2$, although some of the derivations carry over to the general $d$-dimensional case.
We consider the specific case of Example~\ref{expower}, where the functions in Liebscher's construction are power functions, $g_j^{(k)}(t)=t^{p_j^{(k,K)}}$ with $p_j^{(k,K)}\in[0,1]$, $j\in\{1,2\}$ as in~\eqref{eq:p-from-a}. Assuming that $K$ is fixed and limiting ourselves to $d=2$, we denote for notation simplicity 
$p_k\coloneqq p_1^{(k,K)}$ and $q_k\coloneqq p_2^{(k,K)}$ for $k\in\{1,\dots,K\}$.
Recall that, in view of~\eqref{eq:lieb_assumption}, $\sum_{k=1}^K p_k=\sum_{k=1}^Kq_k=1$.
Under the above assumptions, the \LF copula denoted by $\CLF $ has the form 

\begin{equation}\label{eq:Liebscher_Frechet}
    \CLF (u,v)=\prod_{k=1}^K \min(u^{p_k},v^{q_k}), \quad\quad (u,v)\in[0,1]^2.
\end{equation}

In the particular case where $K=2$, it is referred to as the BC2 copula by~\cite{mai2011bivariate} and
it is proved that any 
bivariate extreme-value copula with arbitrary discrete dependence measure can be represented as the geometric mean of 
BC2 copulas, which corresponds to the situation where $K$ is even in~(\ref{eq:Liebscher_Frechet}).
We also refer to~\cite{trutschnig2016mass} for further links with extreme-value theory.

\subsection{Geometric description of $\CLF$}

\begin{figure}
    \centering
    	\input{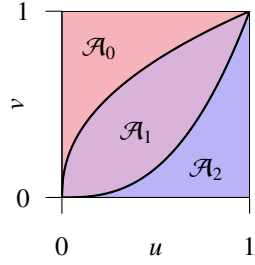}
    \caption{Partition of the unit square defined in \eqref{eq:partition}, for the copula of Equation~\eqref{eq:LF-example}. }
    \label{fig:LF-example}
\end{figure}

For all $k\in\{1,\dots,K\}$, introduce $r_k=p_k/q_k\in[0,\infty]$. For notation simplicity, we shall let $r_0=0$ and $r_{K+1}=\infty$.
Since the above product~(\ref{eq:Liebscher_Frechet}) is commutative, one can assume without loss of generality that
the sequence $(r_k)_{0\leq k\leq K+1}$ is nondecreasing. The copula $\CLF $ can be easily expressed on the partition of the unit square
$[0,1]^2$ defined by the following moon shaped subsets (see the illustration in  Fig.~\ref{fig:LF-example} with $K=2$)
\begin{equation}\label{eq:partition}
 \mathcal{A}_k=\left\{(u,v)\in[0,1]^2 : u^{r_{k+1}}<v\leq u^{r_k}\right\}, \quad\,\,k\in\{0,\ldots,K\}.   
\end{equation}
\begin{Prop}\label{prop:liebFre}
Let $\CLF$ be the \LF copula defined in~\eqref{eq:Liebscher_Frechet}. Then, for all $(u,v)\in[0,1]^2$,
\begin{equation}\label{eq:LF_epxression}
    \CLF (u,v)= \sum_{k=0}^K u^{1-\bar{p}_k}v^{\bar{q}_k}\mathbb{1}[(u,v)\in \mathcal{A}_k],
\end{equation}
where $\bar{x}_k = x_1+\ldots+x_k$, with the convention that $\bar{x}_0=0$. 
Moreover, the singular component of $\CLF $ is
\begin{equation}\label{eq:LF_singular}
    \SLF (u,v)=\sum_{k=1}^K\min(p_k, q_k)\min(u,v^{1/r_k})^{\max(1, r_k)}.
\end{equation}
The singular component $\SLF$ and the absolute continuous component $\ALF =\CLF -\SLF $ weights are   $\sum_{k=1}^K \min(p_k, q_k)$ and $1-\sum_{k=1}^K \min(p_k, q_k)$, respectively.
\end{Prop}

\noindent A key property of the \LF copula~\eqref{eq:Liebscher_Frechet} is the presence of multiple singular components
lying on the curves $v=u^{r_k}$ with associated weights $\min(p_k, q_k)$, $k\in\{1,\dots,K\}$.  As an illustrative example, let us consider the bivariate \LF copula defined with $p_1 = 1-p_2 = 1/3$, $q_1=1-q_2=3/4$,
\begin{equation}\label{eq:LF-example}
    \CLF (u,v)=\min(u^{1/3},v^{3/4})\min(u^{2/3},v^{1/4}), \quad\quad (u,v)\in[0,1]^2,
\end{equation}
which entails $r_0=0, r_1=4/9, r_2=8/3$ and $r_3=\infty$. The moon shaped subsets of the partition of the unit square
$[0,1]^2$ are represented on Fig.~\ref{fig:LF-example}, 
and the expressions of $\CLF$ and the singular component $\SLF$ are as follows:
\begin{align*}
    \CLF (u,v)=\left\{
    \begin{array}{ll}
        u &\text{on } \mathcal{A}_0,\\
        u^{2/3}v^{3/4} &\text{on } \mathcal{A}_1,\\
        v &\text{on } \mathcal{A}_2,\\
    \end{array}\quad\quad 
    \right.\text{and}\quad\quad 
    \SLF(u,v) = \frac{1}{3}\min\big(u,v^{9/4}\big)+
                \frac{1}{4}\min\big(u^{8/3},v\big).
\end{align*}

\begin{figure}[ht]
\centering
  \includegraphics[width=0.29\textwidth]{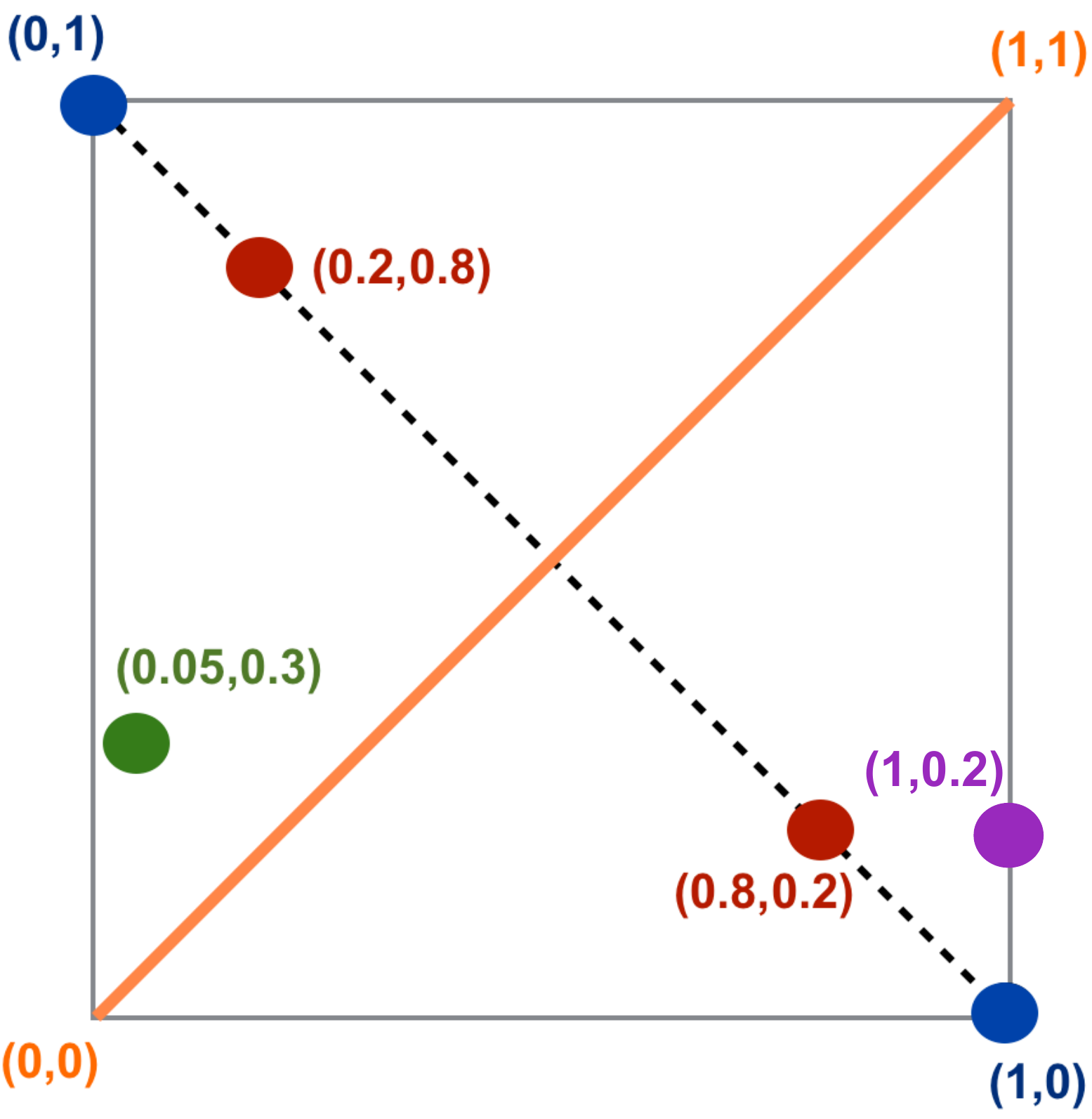}
  \includegraphics[width=0.3\textwidth]{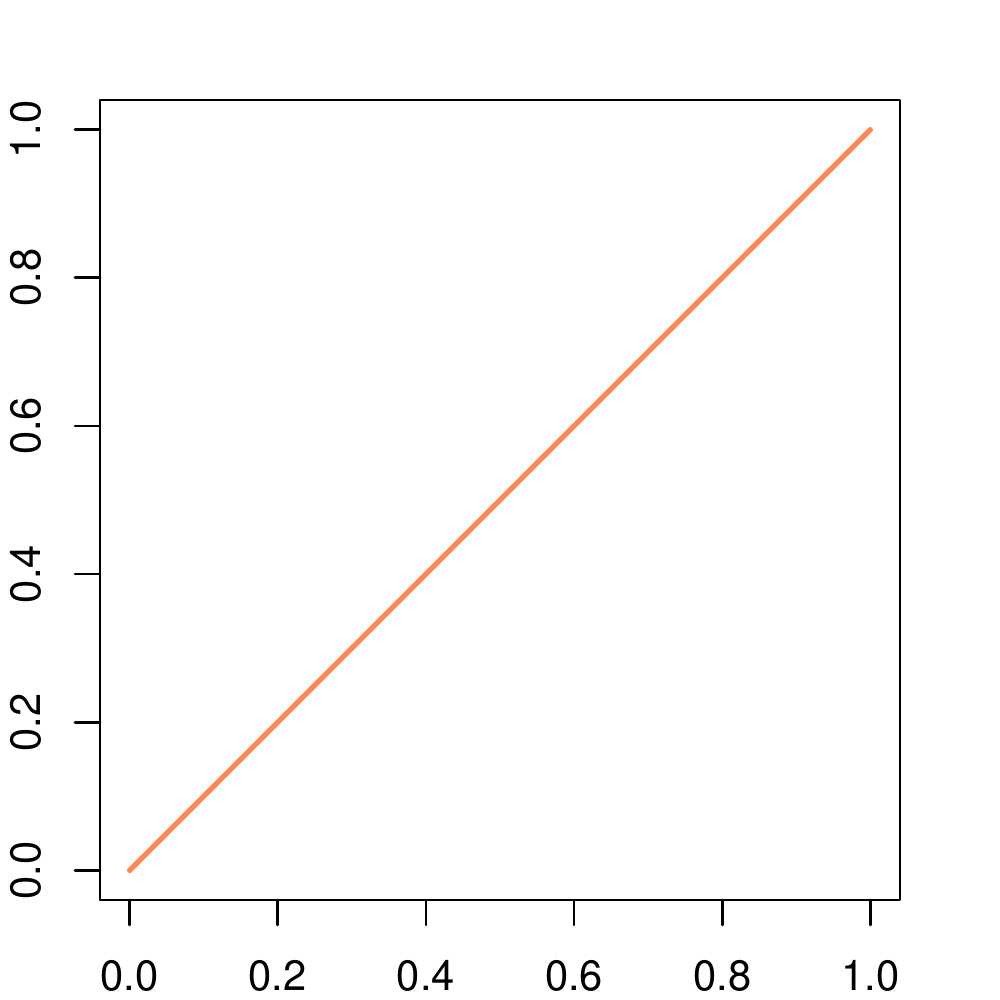}
  \includegraphics[width=0.3\textwidth]{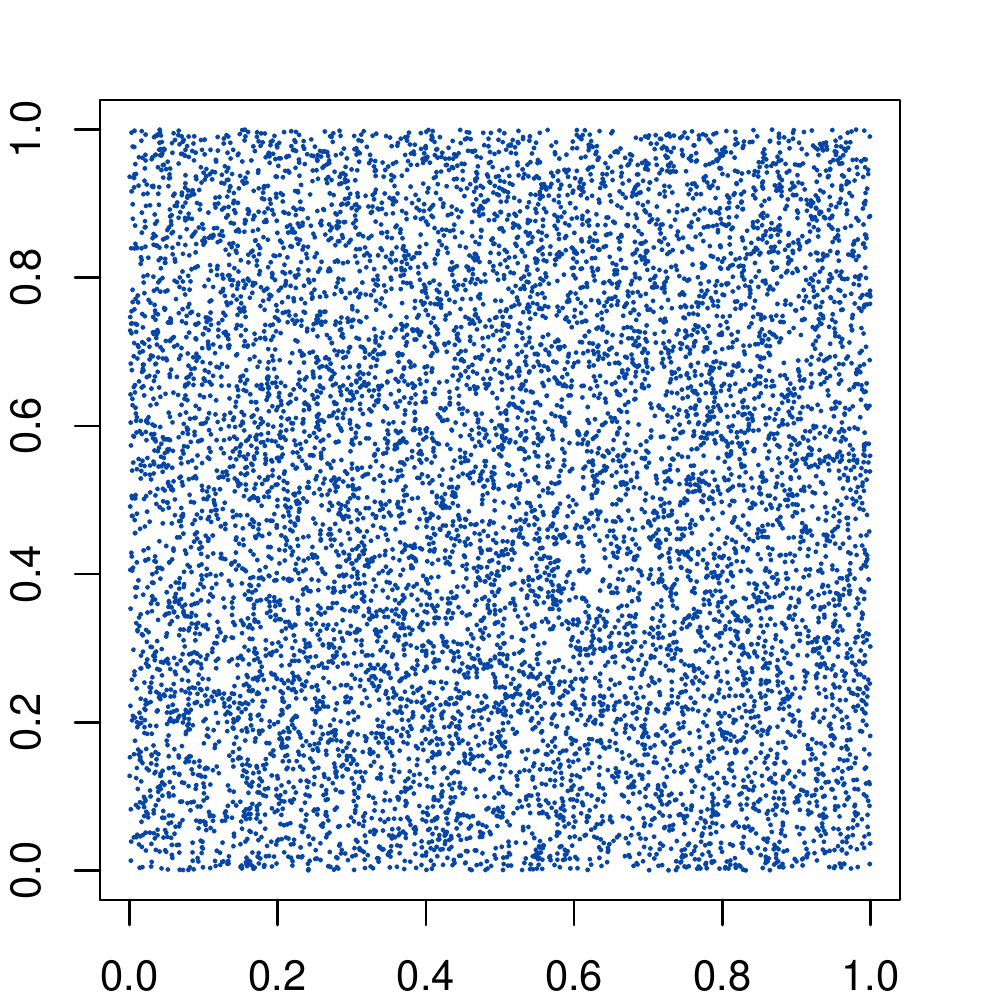}
  \includegraphics[width=0.3\textwidth]{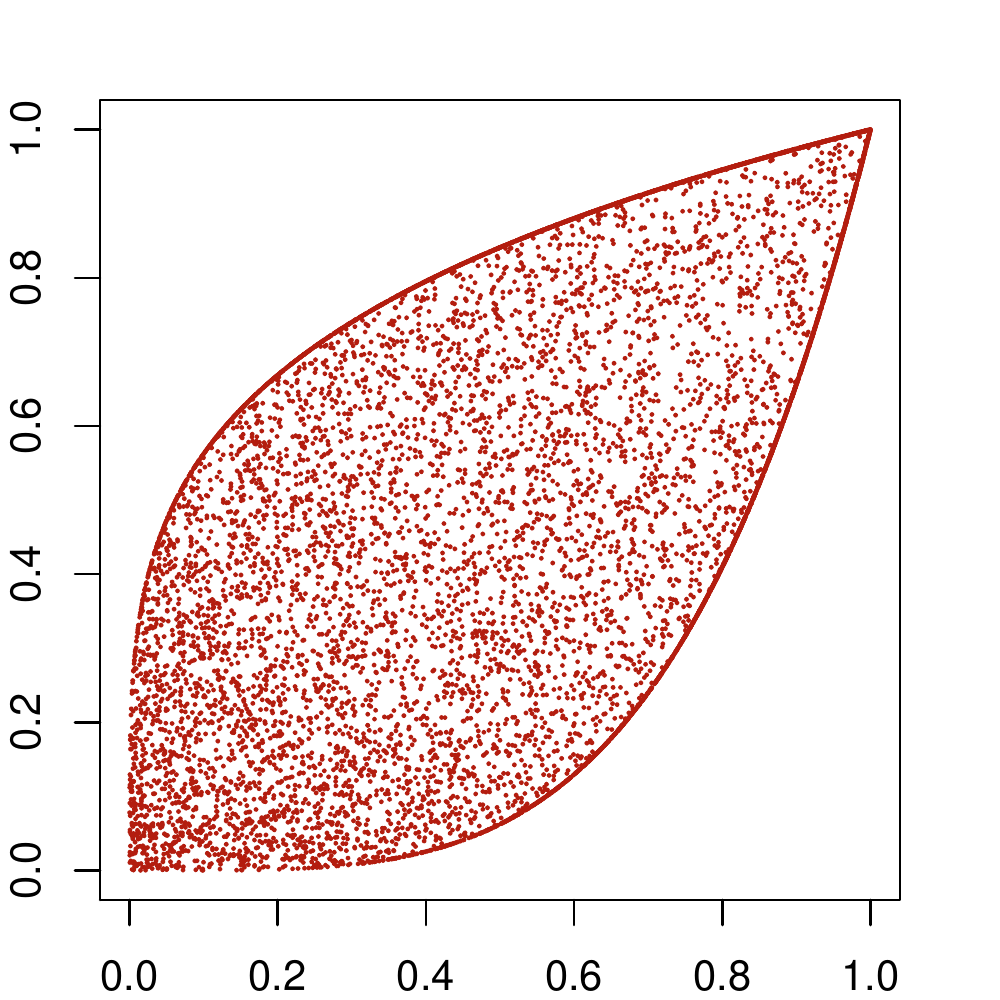}
  \includegraphics[width=0.3\textwidth]{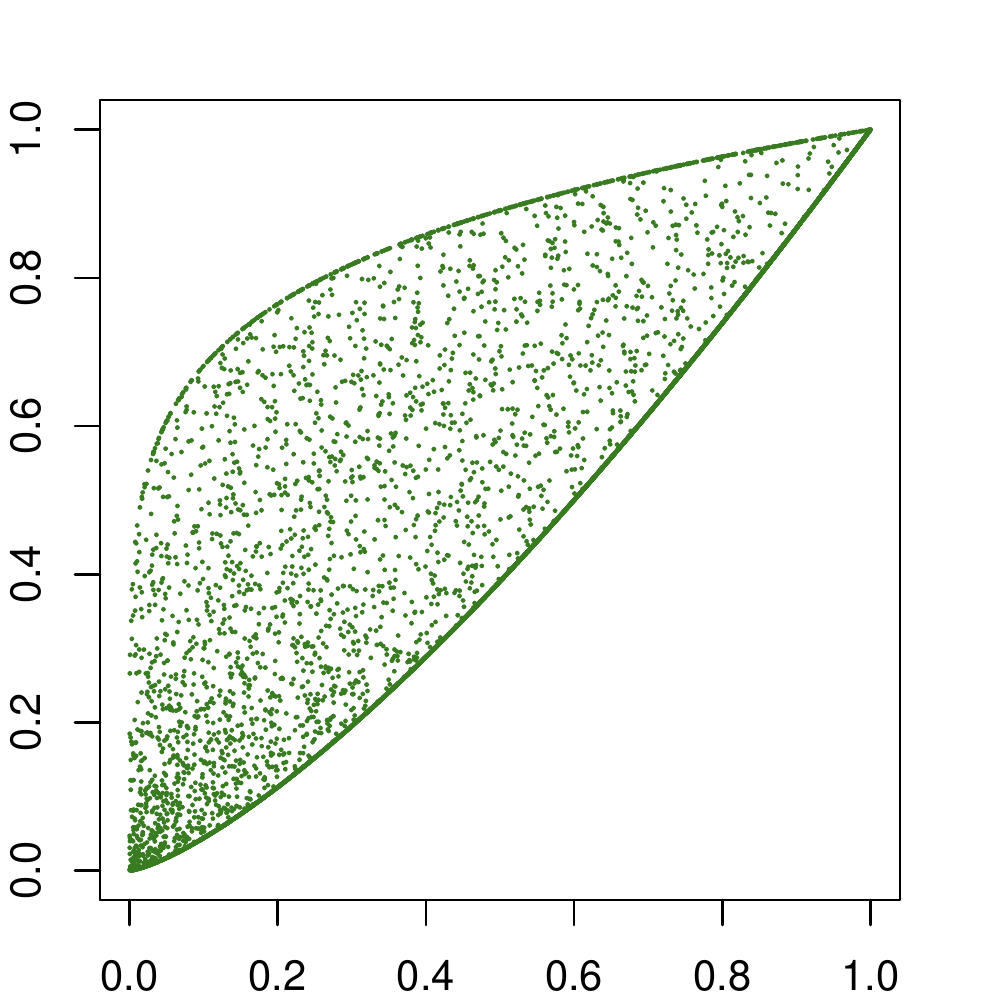}
  \includegraphics[width=0.3\textwidth]{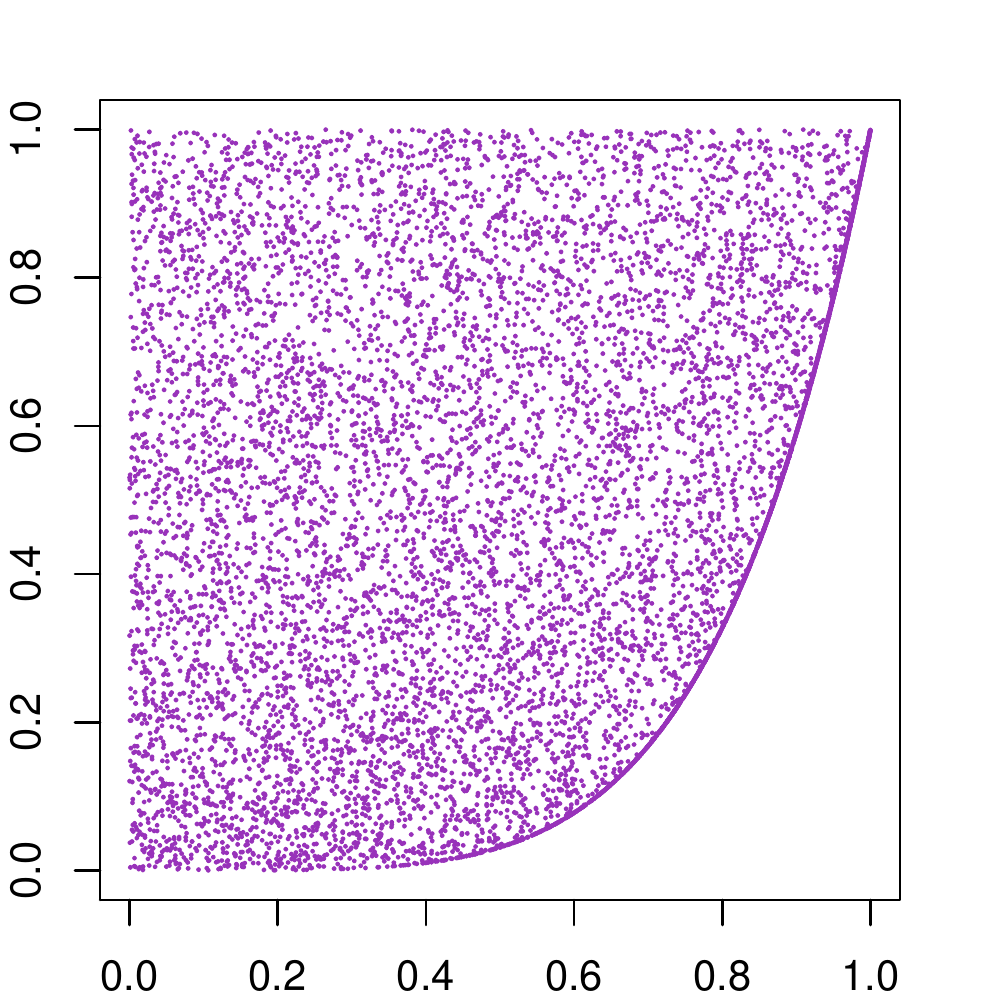}
 \caption{Top-left: representation of the $p\times q$ square unit space. Other five panels: scatter plots of $n=10^4$ data points sampled from \LF copula with $K = 2$. Choices for parameters $(p,q)$ (such that $p_1=p$,  $p_2=1-p$,  $q_1=q$,  $q_2=1-q$) are summarized on the top-left panel. Complete dependence (top-middle),  complete independence (top-right), symmetric (bottom-left), asymmetric (bottom-middle), degenerate asymmetric (bottom-right).
 }
 \label{fig:cases1_2sampled}
 \end{figure}

\begin{figure}[ht]
 \centering
  \includegraphics[width=0.3\textwidth]{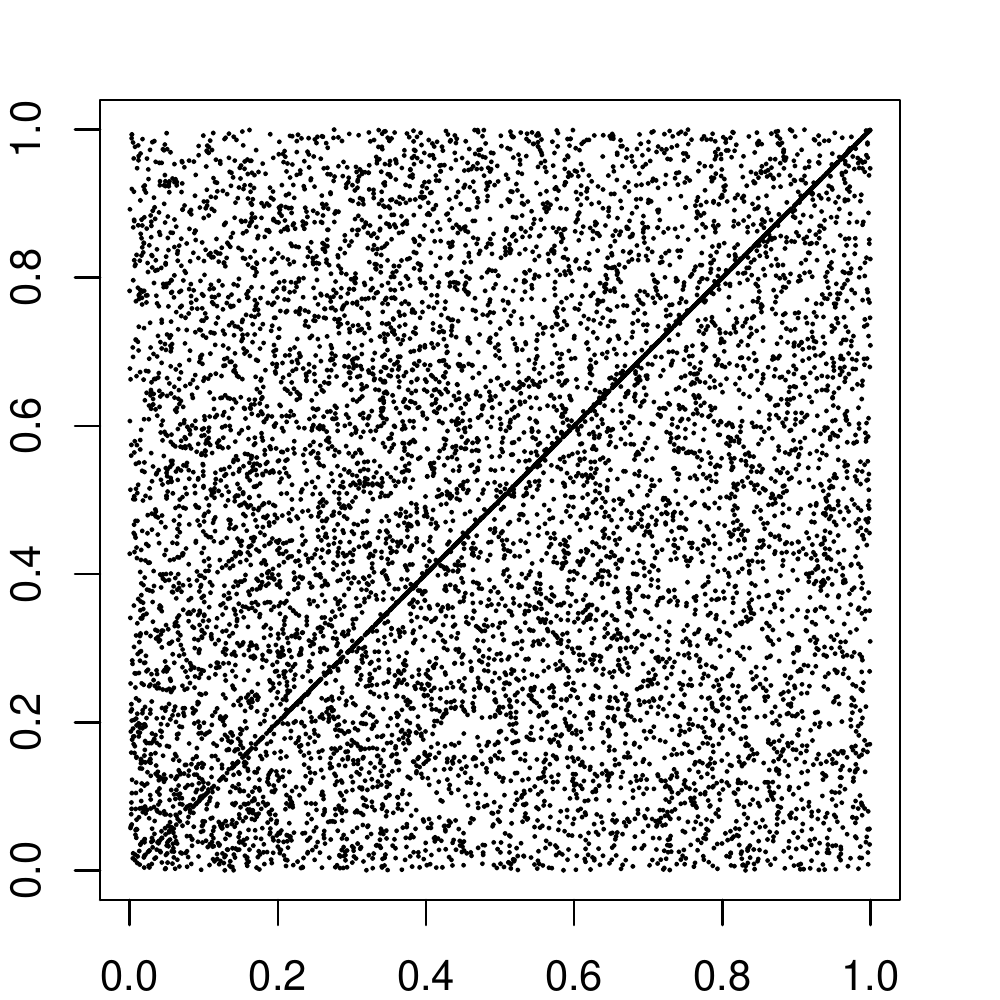}
  \includegraphics[width=0.3\textwidth]{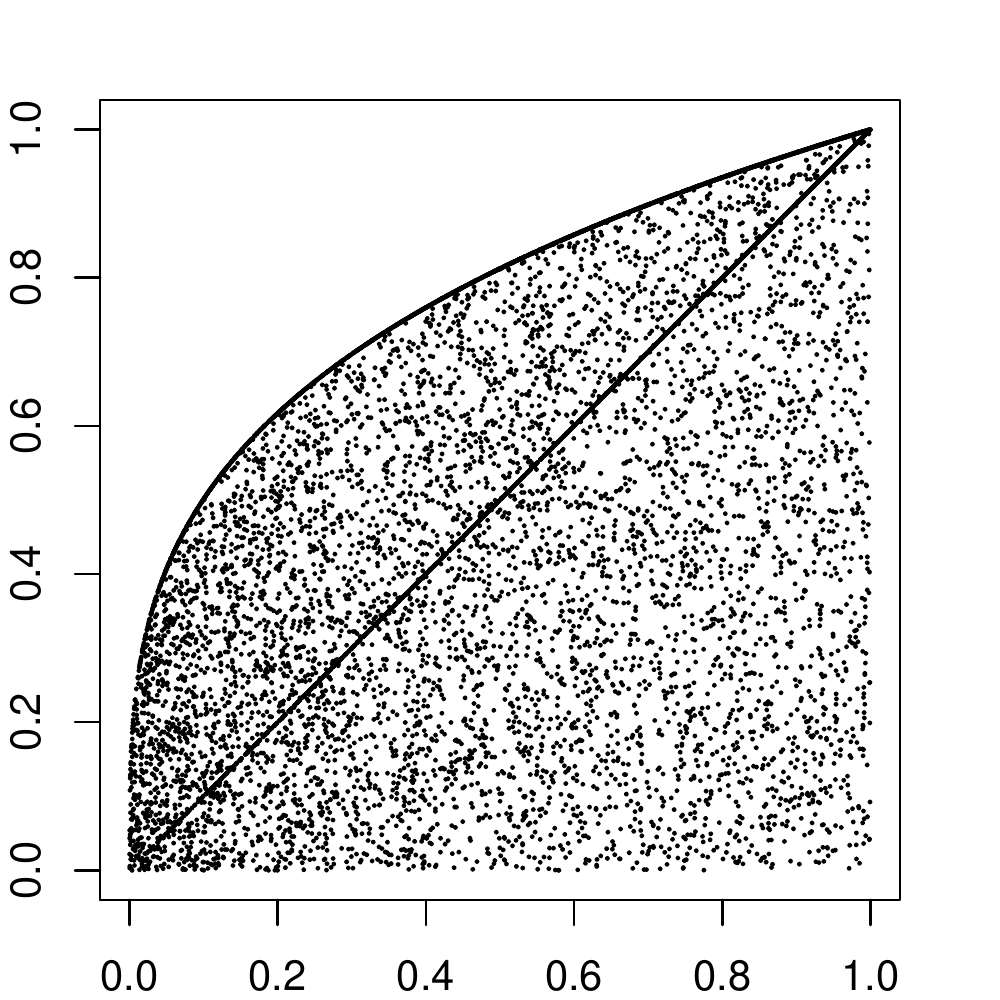}
  \includegraphics[width=0.3\textwidth]{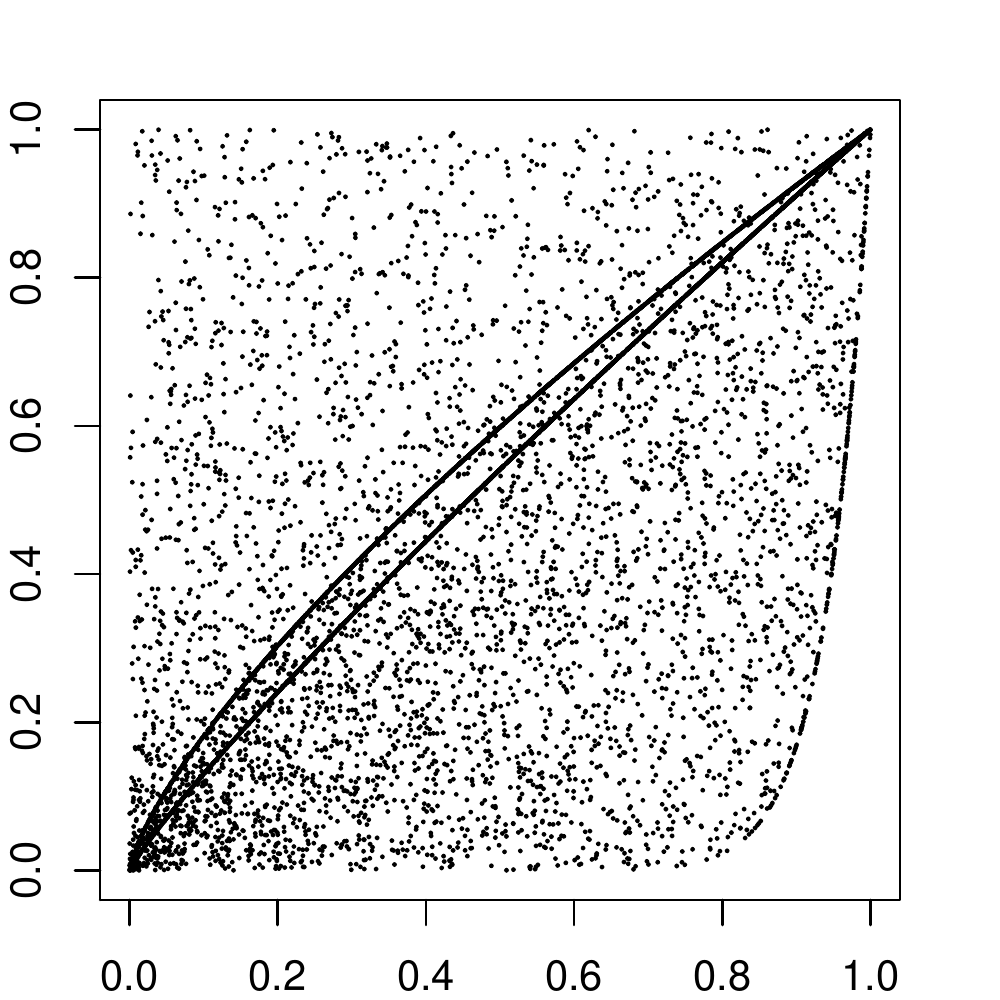}
  \includegraphics[width=0.3\textwidth]{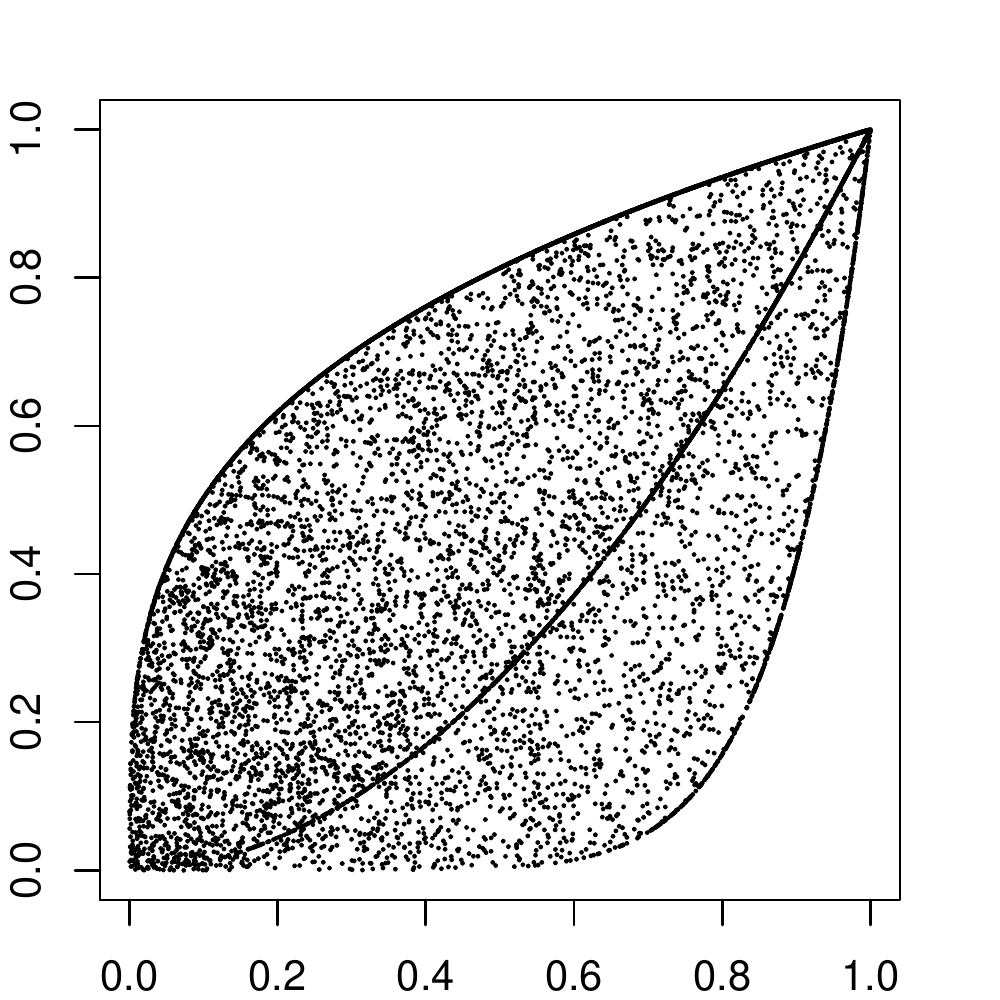}
 \includegraphics[width=0.3\textwidth]{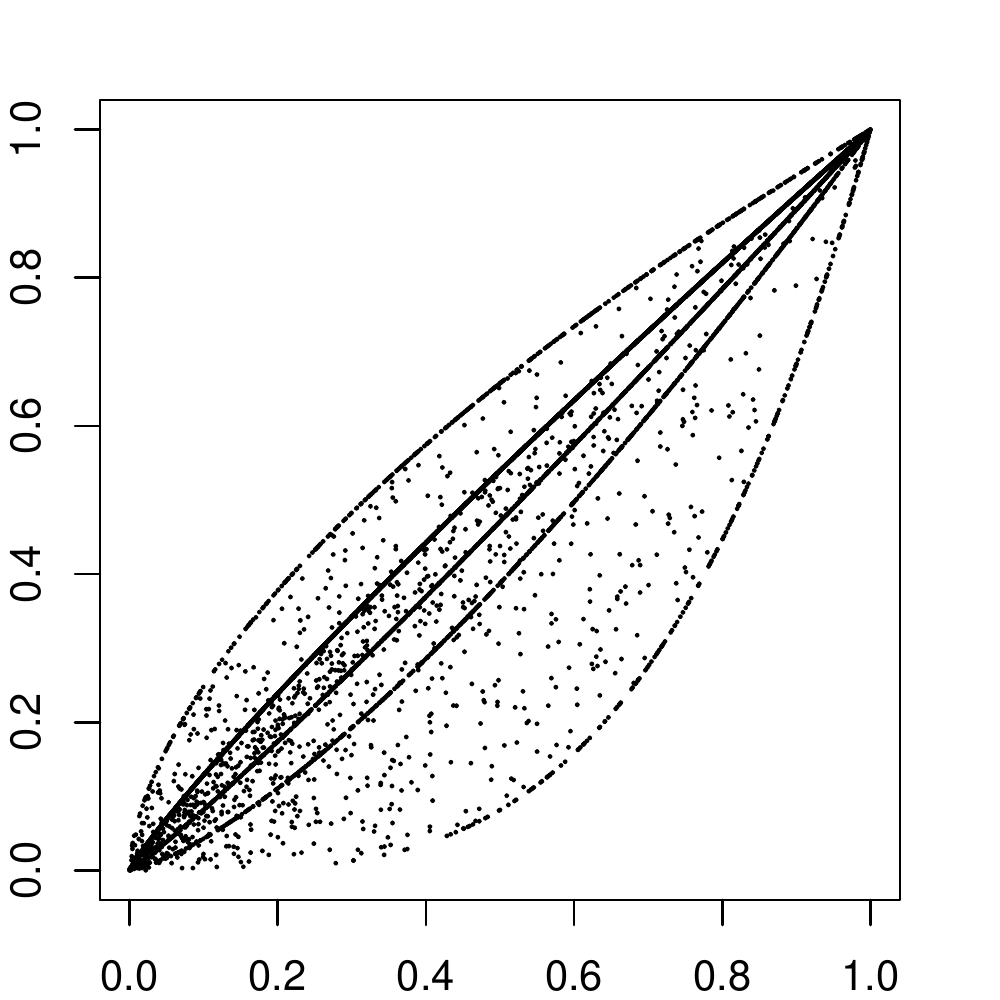}
  \includegraphics[width=0.3\textwidth]{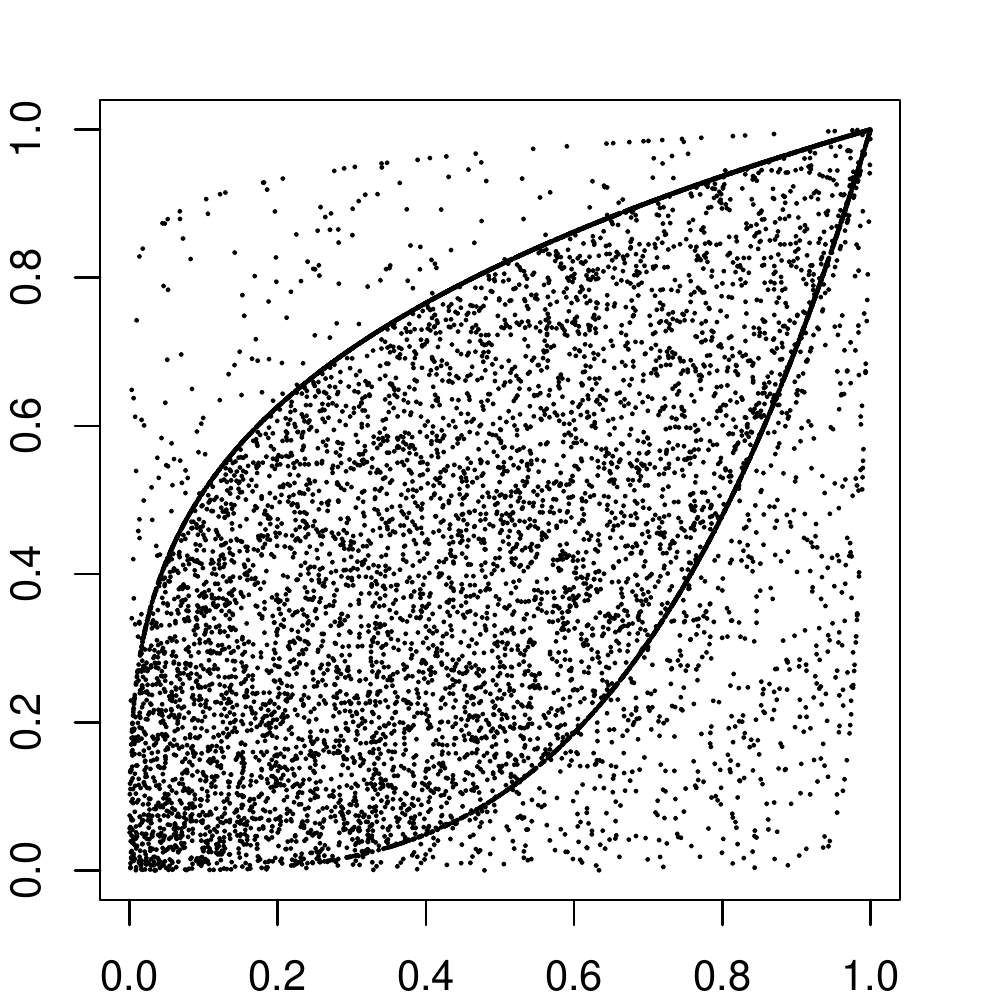}
 \caption{Scatter plots of $n=10^4$ data points sampled from \LF copula with $K>2$. Top-left: $K=3$, $r_k\in\{0, 1, \infty\}$; Top-middle: $K=3$, $r_k\in\{0.3, 1,  \infty\}$; Top-right: $K=4$, $r_k\in\{0,  0.7,  0.9, 17\}$; Bottom-left: $K=4$, $r_k\in\{0.3,  1.9,  8.3, 8.3\}$; Bottom-middle: $K=5$, $r_k\in\{0.6, 0.9, 1.1, 1.4, 3.6\}$; Bottom-right: $K=6$, $r_k\in\{0.04,  0.3,  2.8,  3.3,  4.2, 58\}$.}
 \label{fig:cases1_2sampled2}
 \end{figure}

Appropriate choices of pairs $(p_k,q_k)_{k=1,\dots,K}$
may lead to a number of singular components ranging from $0$ to $K$.
The independence copula (no singular component) is obtained for instance with $p_i=b_j=1$ for a given pair $(i,j)$, $i\neq j$, 
the comonotonic copula (one singular component) is obtained by choosing $p_k=q_k,\,\forall k\in\{1,\ldots,K\}$,
and a copula with exactly $K$ singular components can be obtained provided that $0<r_1<r_2<\dots<r_K<\infty$. See Fig.~\ref{fig:cases1_2sampled} and Fig.~\ref{fig:cases1_2sampled2} for illustrations. 
As a comparison, \citet{cuadras1981continuous} copula given by $C(u,v)=(uv)^{1-\theta}\min(u, v)^\theta$, $\theta\in[0,1]$ is limited to
a single singular component, necessarily on the diagonal $v=u$.
Similarly, \citet{marshall1967generalized} copula is defined by $C(u,v) = \min(u^{1-\alpha}v, uv^{1-\beta})$, $(\alpha,\beta)\in[0,1]^2$
and has only one singular component located on the curve $v=u^{\alpha/\beta}$.
The proposal by~\citet{lauterbach2015some} based on singular mixture copulas includes \LF copula~\eqref{eq:Liebscher_Frechet} in the particular case
when $K=2$ but is limited to two singular components.
Finally, Sibuya copulas~\citep{hofert2013sibuya} is a very general family of copulas: Let us point out that, in the bivariate case, a non-homogeneous Poisson Sibuya copula  allows for only one singular component, this singular component being supported 
by a curve with very flexible shape (see Remark~4.2 in the previously referenced work for further details).

\subsection{Dependence and association properties of $\CLF$}
\label{subsec:FreProp}

We consider here several measures of dependence and measures of association between the components of the bivariate \LF copula~\eqref{eq:Liebscher_Frechet}. Some of these measures are already dealt with  in great generality in Section~\ref{sec:properties}, while some others seem to be   tractable only in the \LF copula  case: 
 Blomqvist's medial correlation coefficient, 
Kendall's $\tau$ 
and Spearman's $\rho$.

\paragraph{Tail dependence}
Recall that for the comonotonic copula $C_{\text{C}}$, it holds $\Lambda_U(C_{\text{C}};\cdot,\cdot)=\min(\cdot,\cdot)$. Then,
Corollary~\ref{cor:lambda} yields 
\begin{equation*}
\lambda_L(\CLF )= \prod_{k=1}^K \mathbb{1}(p_k=q_k), \quad\quad
\lambda_U(\CLF )= \sum_{k=1}^K \min(p_k,q_k).
\end{equation*}
In other words, the lower tail dependence coefficient is non zero only in the case when $p_k=q_k$ for all $k\in \{1,\ldots,K\}$, where $\CLF$ boils down to the comonotonic copula, while the upper tail dependence coefficient coincides with the weight of the singular component $\SLF$ (see Proposition~\ref{prop:liebFre}).
 These results were also established in~\cite[][Lemma~3]{mai2011bivariate}, in the particular case where $K=2$.

\paragraph{Dependence}

It is well-known that the comonotonic copula $\CF$ fulfills the following positive dependence properties defined in Section~\ref{subsec:dependence}, namely it is \TPtwo, \PQD, \LTD and \SI \citep{Nelsen}. According to Proposition~\ref{prop:long_list_of_dependence}, we can thus conclude that the \LF copula $\CLF$ also satisfies these  positive dependence properties.

\paragraph{Stability properties} It is easily seen that \LF copula~\eqref{eq:Liebscher_Frechet} is max-stable. Proposition~\ref{propstable2} thus entails that this copula is stable with respect to Liebscher's construction. Another consequence is that \LF construction~(\ref{eq:Liebscher_Frechet}) can be interpreted as a possible cdf for modelling bivariate maxima.

\paragraph{Dependence coefficients}
 The $\beta$-Blomqvist's medial correlation coefficient (\citep{Nelsen}, Paragraph~5.1.4) defined by $\beta(C) = 4C\big(\frac{1}{2},\frac{1}{2}\big)-1$,  
Kendall's $\tau$ (\citep{Nelsen}, Paragraph~5.1.1) 
and Spearman's $\rho$ (\citep{Nelsen}, Paragraph~5.1.2) 
 defined by 
$$
\tau(C) =  4\int_{[0,1]^2}  C(u,v)\,\ddr C(u,v)-1 \;\mbox{ and }\;
 \rho(C) = 12\int_{[0,1]^2}  C(u,v)\,\ddr u\ddr v-3 
$$
are provided in next proposition.
\begin{Prop}
\label{prop3coeffs}
Blomqvist's medial correlation coefficient, 
Kendall's $\tau$ 
and Spearman's $\rho$ for the \LF copula~\eqref{eq:Liebscher_Frechet} are respectively given by 
\begin{align*}
    \beta(\CLF ) &=  2^{\sum_{k=1}^K \min(p_k, q_k)}-1,\\
    \tau(\CLF ) &= 1 - \sum_{k=1}^{K-1} \frac{(1-\bar p_k) \bar q_k (r_{k+1}-r_k)}{(\bar q_k r_k + (1-\bar p_k))(\bar q_k r_{k+1} + (1-\bar p_k))}, \\
    \rho(\CLF ) &= \frac{12(1+r_1+r_1 r_K)}{(2+r_1)(1+2r_K)} -3 + 
    \sum_{k=1}^{K-1} \frac{r_{k+1}-r_k}{((1+\bar q_k) r_k + (2-\bar p_k))((1+\bar q_k) r_{k+1} + (2-\bar p_k))},
\end{align*}
where $\bar{x}_k = x_1+\ldots+x_k$.
\end{Prop}
\noindent It appears that $\beta$-Blomqvist medial correlation coefficient is closely related to the upper tail dependence coefficient:
$\beta(\CLF )=2^{\lambda_U(\CLF )}-1$. Besides, in the particular case where $K=2$,  these results coincide the ones
of Lemma~2 of~\cite{mai2011bivariate}: the Kendall's $\tau$ can be simplified as $\tau(\CLF )=p_1+q_2=\lambda_U(\CLF )$.
No similar simplification seems to be possible for Spearman's $\rho$.

For the special case of copula~\eqref{eq:LF-example}, we have $\lambda_L(\CLF) = 0$, 
 $\lambda_U(\CLF) = \tau(\CLF) = p_1+q_2 = \frac{1}{3}+\frac{1}{4} = \frac{7}{12}\approx 0.583$, 
$\beta(\CLF)=2^{7/12}-1\approx 0.498$ and $\rho(\CLF)\approx 0.298$.

\subsection{Iterative construction for $\CLF$}

Algorithm~\ref{algo:Iterative_copula_sampling} can be simplified when specified to \LF setting since: (i) sampling from the comonotonic copula is straightforward and only requires sampling from the uniform distribution $\mathcal{U}(0,1)$, 
and (ii) power functions benefit from an explicit inverse. The specific sampling procedure for this construction is described in detail as Algorithm~\ref{algo:sampleFrechet}. 

\vspace{0.5cm}
\begin{center}
\begin{minipage}{10.2cm}
\begin{algorithm}[H]
\begin{footnotesize}
\DontPrintSemicolon
\textbf{Input} $\left[a_j^{(k)}\right]_{(k,j)}$
\algocomment{exponents of power functions $f_j^{(k)}$}\;
$X^{(1)}_j \sim \mathcal{U}(0,1)$
for each $j=1,\dots,d$\;
\For{$k=2,\ldots,K$}{
\textbf{Sample} $Y\sim \mathcal{U}(0,1)$, independently of $X^{(k-1)}_1,\ldots,X^{(k-1)}_d$\;
\For{$j\in\{1,\ldots,d\}$}{
 \textbf{Compute} $X^{(k)}_j = \text{max}\left(\left(X_j^{(k-1)}\right)^{\frac{1}{1-a_j^{(k)}}}, Y^{\frac{1}{a_j^{(k)}}}\right)$
 }}
 \textbf{Output} $\bm{X}=(X_1^{(K)},\ldots,X_d^{(K)})\sim \CLF $
\;
\caption{Iterative sampling scheme for \LF copula~\eqref{eq:Liebscher_Frechet} }\label{algo:sampleFrechet}
\end{footnotesize}
\end{algorithm}
\end{minipage}
\end{center}


\section{Bayesian inference}\label{sec:simu}
In this section, we provide a simple strategy to make Bayesian inference on any \LC copula based on an Approximate Bayesian computation algorithm (ABC, see for instance \citep{abc,robert2018abc,karabatsos2018approximate} for reviews).
ABC is a ``likelihood-free'' method usually employed for inference of models with intractable likelihood: it enables to perform approximate Bayesian analysis on any statistical model from which it is possible to sample new data, without the need to explicitly evaluate the likelihood function.

Let   $\bm{X}_\obs=\{\bm X_{\obs,1},\ldots,\bm X_{\obs,n}\}$ be the observed data, where $\bm X_{\obs,i}=(X_{\obs,i,1},\ldots,X_{\obs,i,d})$, $i\in\{1,\ldots,n\}$, and assume that the statistical model for $\bm{X}_\obs$ is described by a likelihood function $\mathcal{L}_\theta$ with parameter $\theta$ which is to be inferred. The basic scheme of one step of ABC is the following: 
\begin{enumerate}
\item Sample $\theta$ from the prior distribution $\pi(\theta)$;
\item Given $\theta$, sample $\bm{X}_1,\ldots,\bm{X}_n$ from $\mathcal{L}_{\theta}$, and set $\bm{X}=\{\bm X_1,\ldots,\bm X_n\}$; 
\item If  $\bm{X}$ is too different from $\bm{X}_\obs$, discard $\theta$, otherwise, keep $\theta$. 
\end{enumerate}
The outcome of the ABC algorithm is a sample of values of the parameter $\theta$ approximately distributed according to its posterior distribution. 
The basic (rejection) ABC approach in point 3. amounts to \emph{a priori} specifying a tolerance level $\epsilon>0$, and then  keeping $\theta$ if $d(\bm{X},\bm{X}_\obs)<\epsilon$ for some distance $d(\,\cdot\,,\,\cdot\,)$ between samples. 
Another common approach employed in this paper consists in selecting the tolerance level $\epsilon$ as a fixed quantile  of the distances $d(\bm{X},\bm{X}_\obs)$. More specifically, Steps 1. to 3. are repeated $M'$ times, out of which $M$ are retained, yielding a quantile of order $M/M'$ \citep{robert2018abc}. In other words, the $M$ retained parameters are those associated with the smallest values of the distance $d(\bm{X},\bm{X}_\obs)$.
In this paper, we choose as distance between samples the Hilbert distance introduced by \cite{2017arXiv170105146B}, henceforth denoted by $d_H(\,\cdot\,,\,\cdot\,)$. The Hilbert distance is an approximation of the Wasserstein distance between empirical probability distributions which preserves the desirable properties of the latter in the context of ABC, while being considerably faster to compute in multivariate data settings. More precisely, given two samples, $y_{1:n}$ and $z_{1:n}$, the Hilbert distance of order 1 associated with the Euclidean distance (henceforth simply referred to as the Hilbert distance) is defined as $d_H(y_{1:n},z_{1:n})=\frac{1}{n}\sum_{i=1}^{n} ||y_i-z_{\sigma(i)}||$, where $\sigma(i)= \sigma_z \circ \sigma^{-1}_y(i)$ for all $i=1,\dots,n$, and  $\sigma_y$ and $\sigma_z$ are the permutations obtained by mapping the
vectors $y_{1:n}$ and $z_{1:n}$ through their projection via the Hilbert space-filling curve \citep{sagan} and sorting the resulting vectors in increasing order (see \citep{2017arXiv170105146B} for details).

The choice for ABC is motivated by two main reasons. First, it is nontrivial in general to derive the likelihood of copulas, especially for Liebscher copulas which involve differentiating a product of $K$ terms. All the more, the specific case of \LF copulas induces up to $K$ singular components which precludes  a general evaluation of the likelihood. 
Second, sampling new data  (step 2. above) from a Liebscher copula is straightforward and fast thanks to the iterative procedure of Algorithm \ref{algo:Iterative_copula_sampling} (Section \ref{sub:constr}). 

Section \ref{lab:priors} introduces the ABC procedure in the case of the \LC copula~\eqref{eq:lieb_cop} and describes the prior distributions on the model parameters. 
The methodology is then illustrated on two data generating distributions: Section~\ref{sub:Inf_liebFre} focuses on the well-specified setting where the data are sampled from the \LF copula; we show that the estimation procedure performs well in this case. Then, Section~\ref{sub:Inf_liebFreNoi} investigates the misspecified setting where the data are sampled from a noisy version of the model; we show that the estimation procedure still performs well, but the estimation accuracy may deteriorate for too large values of the noise. Finally, Section~\ref{sub:abc_VS_mle} compares our proposed ABC approach to a likelihood-based technique.

\subsection{ABC inference for \LC copulas}\label{lab:priors}

The description of the ABC procedure is first completed by specifying the prior distributions on the model parameters. For simplicity, we here focus on the case of the $d$-dimensional \LC copula with the functions $f_j^{(k)}(\cdot)$, $k\in\{1,\ldots,K\}$, $j\in\{1,\ldots,d\}$ of Algorithm \ref{algo:Iterative_copula_sampling} chosen as the power functions~\eqref{eq:def_powers} introduced in Example~\ref{expower}. 
The parameters of the $f_j^{(k)}(\cdot)$ functions are collected in a $K\times d$ matrix ${A}=[a_j^{(k)}]_{k,j}$, where ${a}^{(1)}=(a_1^{(1)},\ldots,a_d^{(1)})=(1,\ldots,1)$. 
Since the $(K-1)d$ free parameters are constrained to $a_j^{(k)}\in(0,1)$ for $2\leq k\leq K$, we simply choose, by symmetry, independent and uniform distributions $a_j^{(k)}\simiid \mathcal{U}(0,1)$. More flexible distributions like the Beta distribution could be thought of in order to reflect some prior knowledge on these parameters. Additionally, note that different functions $f_j^{(k)}(\cdot)$ would simply require setting prior distributions adapted to the parameters used.

The number of iterative steps $K$ is also considered as a parameter of the model. 
Independently of parameters $A$, $K$ is assigned a Zipf distribution, $K\sim \text{Zipf}(\xi)+1$, for $\xi>1$. Such a distribution is supported on integers $k\geq 2$ and has probability mass function $\P(K=k)$ proportional to $(k-1)^{-\xi}$. We further choose the parameter $\xi$ to be equal to 2, which insures that $90\%$ of the prior mass for $K$ is supported on most realistic values $2\leq K \leq 6$. This can be changed depending on applications at hand.
Another option, useful in case where some prior information is available on $K$, is to adopt a Binomial distribution (translated, such that $K\geq2$).
The choice of the two hyper-parameters of the Binomial density could then be set as a function of prior knowledge, such as prior mode and confidence, that one may be able to elicit thanks to expert knowledge or previous studies. 

We are now ready to state the main ABC inference procedure as Algorithm \ref{algo:abc}.

{\begin{algorithm}[h!]
\begin{footnotesize}
\DontPrintSemicolon
\textbf{Input} $\bm{X}_\obs, M', M, (C_k)_{k}.$
\; 
\For{$s=1,\ldots,M'$}{
$K^{(s)} \sim \text{Zipf}(\xi)+1$ \algocomment{sample number of iterations in construction~\eqref{eq:lieb_cop}}\;
\For{$j\in\{1,\ldots,d\}$ \text{\emph{and}} $k\in\{2,\ldots,K^{(s)}\}$}{
 $a^{(1)}_{j}=1$\;
 $a^{(k)}_{j} \sim \mathcal{U}(0,1)$ \algocomment{sample copula parameters}\;
}
${A}^{(s)}=[a_{j}^{(k)}]_{j,k}$
\algocomment{set parameters for Liebscher's construction}\;
$\bm X_{1},\ldots,\bm X_{n}\simiid \CL^{(K^{(s)})}$ 
\algocomment{sample data using Algo. \ref{algo:Iterative_copula_sampling} with power functions and ${A}={A}^{(s)}$}\;
$\bm{X}^{(s)}=\{\bm X_{1},\ldots,\bm X_{n}\}$
\algocomment{set synthetic data}\;
$d_H^{(s)}=d_H(\bm{X}_\obs, \bm{X}^{(s)})$
\algocomment{compute Hilbert distances}\;
}
{\textbf{Compute} $d^*$: the quantile of order $M/M'$ of the distances $[d_H^{(s)}]_{s=1}^{M'}$}\;
\textbf{Output} $\{(\bm{X}^{(m)},{A}^{(m)},{K}^{(m)}):\,\,d_H^{(m)}<d^*\}_{m=1}^{M}$
\algocomment{return $M$ parameters with smallest $d_H$ from $\bm{X}_\obs$}\;
\caption{ABC inference for Liebscher copulas}\label{algo:abc}
\end{footnotesize}\end{algorithm}
}
In general, the sequence of copulas $(C_k)_k$  depends on some sequence of parameters $(\gamma_k)_k$. 
In such a case,  Algorithm~\ref{algo:abc} can be easily amended by adding a step consisting in sampling $\gamma_k$ parameters from some prior distribution to be set based on available prior information or expert knowledge. 
In the following section, we focus on the case of \LF copulas, which do not depend on any additional parameter. Thus, the sampling step for new data $\bm X^{(s)}$ in Algorithm~\ref{algo:abc} is performed with the iterative construction of Algorithm~\ref{algo:sampleFrechet} tailored to \LF copulas. 

\subsection{Numerical illustrations with \LF copulas}\label{sec:illustr}
This section provides two illustrations of the inferential procedures described so far. The first investigates a setting where data are sampled from the \LF model, while the second is concerned with observations from a noisy version of it.
The code, implemented in \textsf{R} using the \texttt{copula} package \citep{copulaR} and \texttt{winference} package \citep{winferencecode} for the Hilbert distance implementation \citep{2017arXiv170105146B}, is available at the following link:\\ \href{https://sites.google.com/site/crispinostat/research?authuser=0}{https://sites.google.com/site/crispinostat/research?authuser=0}.

\subsubsection{Well-specified setting: data from \LF copula}\label{sub:Inf_liebFre}
We generate $n=500$ data points from a $2$-dimensional \LF copula~\eqref{eq:Liebscher_Frechet} with $d=2$,  varying values of $K\geq 2$ and of the parameters in the matrix ${A}$, using Algorithm~\ref{algo:sampleFrechet}. The estimation is then performed with the ABC procedure summarized in Algorithm~\ref{algo:abc}.

Our method provides the full (approximate) posterior distribution of the parameters of interest, making possible to select any strategy to summarize them, possibly driven by the application at hand.
One can for instance compute the posterior distribution of the Spearman's $\rho$, and, thanks to the retained samples, any other quantity of interest.
Here, the performance of the estimation procedure is assessed basing on the following three summary statistics:
(i) Kendall's distribution function $\K(t)=\Pr_C({C}(U,V)<t)$, (ii)
Spearman's $\rho$ index of association, introduced in Section~\ref{subsec:FreProp}, and for which an explicit closed form is obtained for \LF copula in Proposition~\ref{prop3coeffs}, and (iii)
an asymmetry measure, since it is a central motivation of the present work. More specifically, the Cram\'er-von Mises test statistics ${\rm E}_C [({C}(U,V)-{C}(V,U))^2]$ defined in \citet{genest2012tests} 
has been selected since it emerged as a powerful statistic to test the symmetry of a copula.
Following the strategy of \citet{genest2012tests}, the approximate p-values associated with the symmetry test, performed both on the observed sample $\bm{X}_\text{obs}$ and on the retained samples $\bm{X}_m, m\in\{1,\ldots,M\}$, are computed on the basis of 250 bootstrap replicates.

 \begin{figure}[h!]
 \centering
 \includegraphics[width=0.4\textwidth]{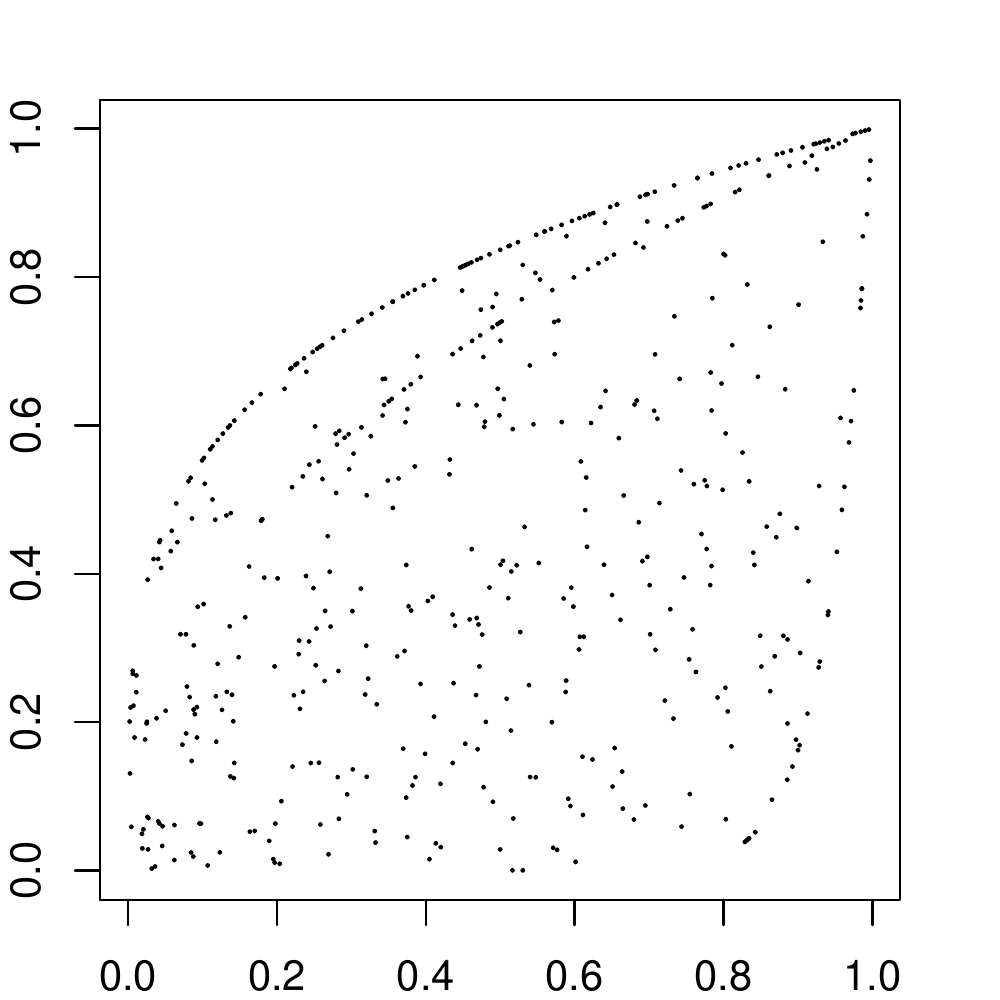}
 \includegraphics[width=0.4\textwidth]{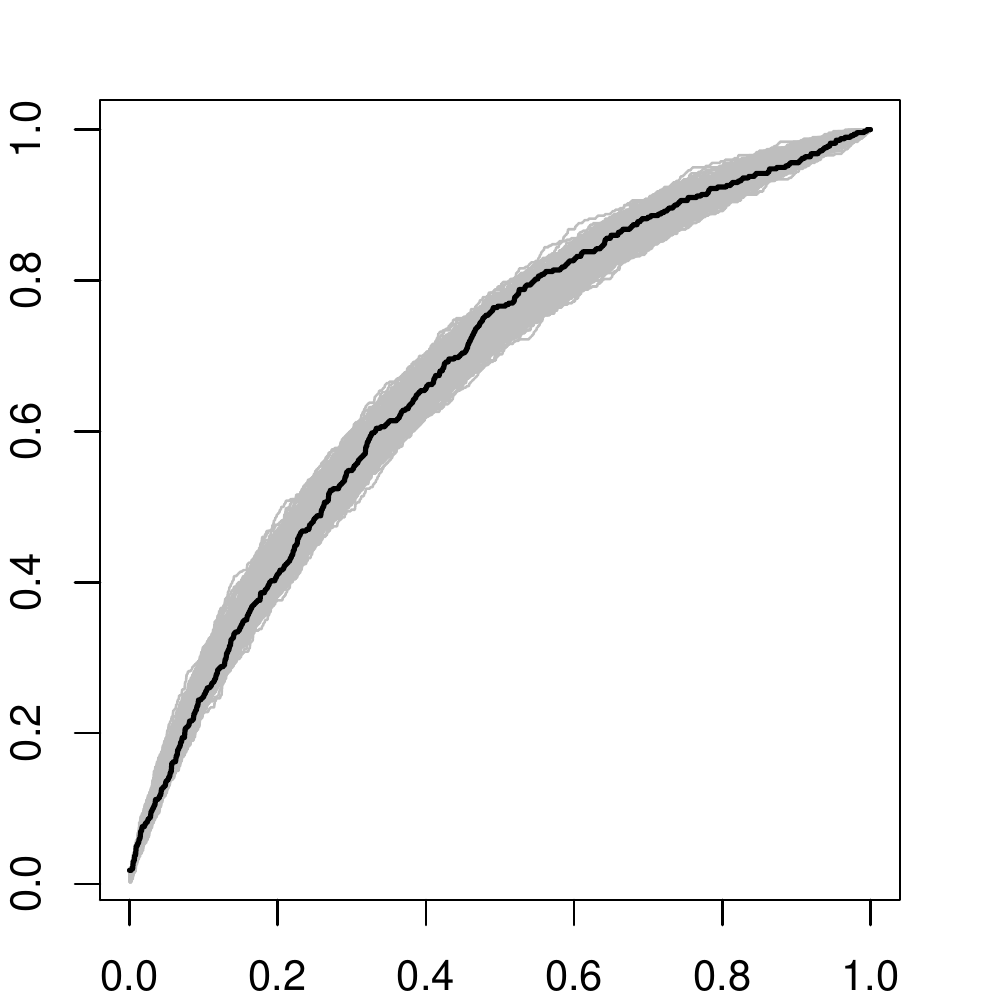}
 \includegraphics[width=0.4\textwidth]{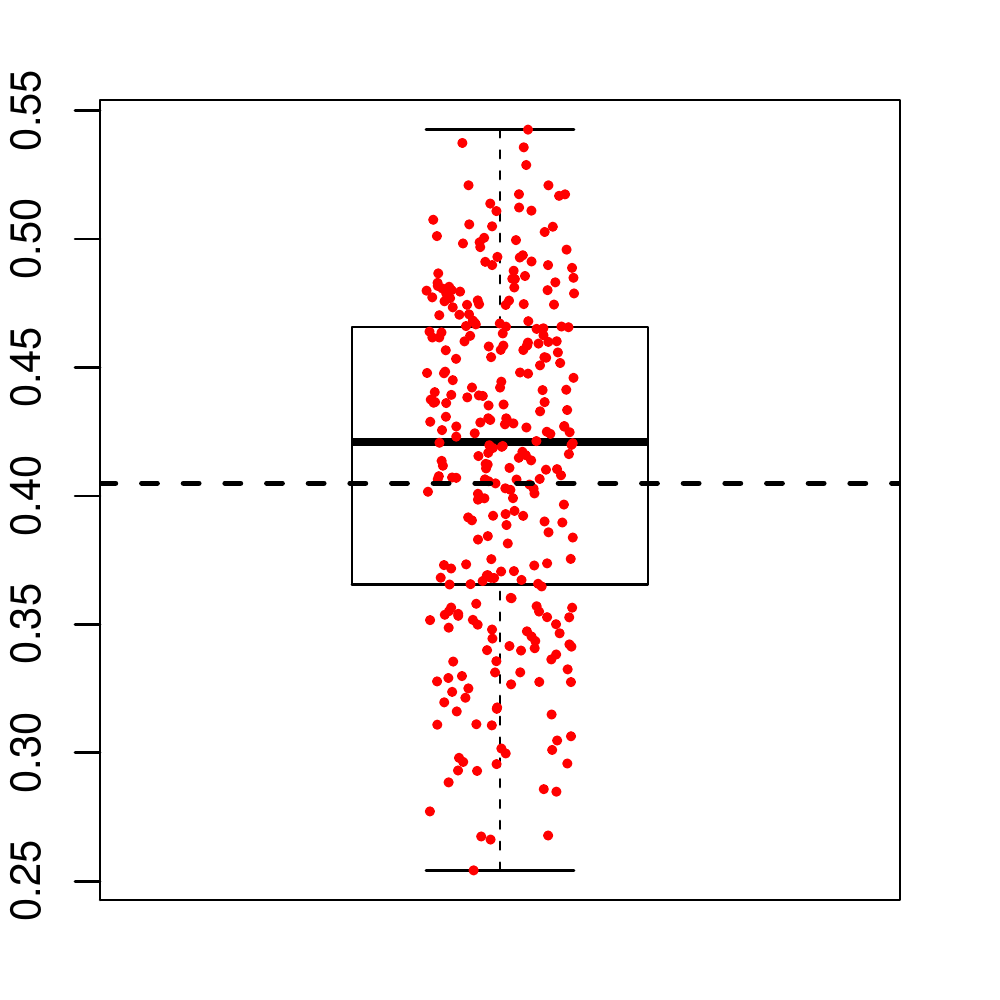}
 \includegraphics[width=0.4\textwidth]{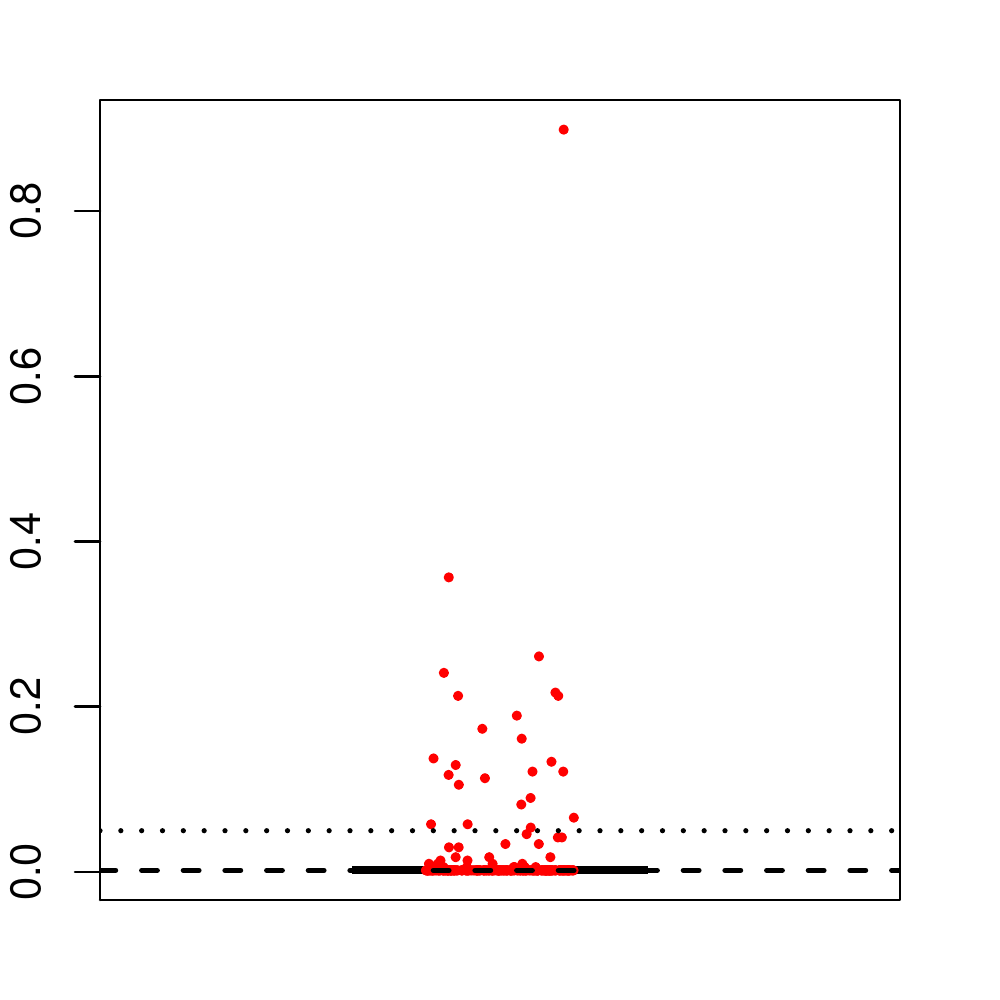}
 \caption{Results from a simulation experiment with $K=3$. 
 Top-left: $n=500$ data points simulated from a \LF copula with $d=2$, $K=3$ and $r_k\in\{0.26, 0.44, 17.32\}$. Top-right: empirical Kendall's distribution function $\hat\K$ of the observed (black) and retained (gray) samples; Bottom-left: boxplot of the posterior distribution of $\rho$ 
 (the dashed line corresponds to the empirical Spearman's $\rho$ of the observed sample); Bottom-right: boxplot of the posterior distribution of the approximate p-values (the dashed line corresponds to $\bm{X}_\text{obs}$, and the dotted line is the 5\% threshold).
 }
 \label{fig:prova}
 \end{figure}

The results obtained on a single simulation experiment are displayed on Fig.~\ref{fig:prova}, where $n = 500$ data points were simulated from the \LF copula~\eqref{eq:Liebscher_Frechet} with $d=2$ and $K=3$ (top-left panel).
The number of ABC iterations was set to $M'=10^4$, of which $M=300$ were retained (resulting in a quantile of order $3\%$). 
The empirical Kendall's distribution functions of the observed  and retained samples are compared 
on the top-right panel; The posterior distribution of $\rho$ is compared to the empirical Spearman's $\rho$ of the observed sample on the bottom-left panel; Finally, the posterior distribution of the approximate p-values is displayed on the bottom-right panel.
Let us highlight that the estimating procedure provides distributions around the true values in the three considered cases. 

We then vary the generating number of iterative steps $K$ and compute the average relative errors $\eta_\K$ and $\eta_\rho$ for $\K$ and $\rho$ between the values computed on the observed sample and on the $M$ samples retained by ABC:
\begin{align}\label{eq:relative-errors}
    \eta_\K=\frac{1}{M}\sum_{m=1}^M \frac{\Vert \hat\K_\obs-\hat\K_m \Vert_1}{\Vert \hat\K_\obs\Vert_1}, \, \quad\text{ and }\quad
    \eta_\rho=\frac{1}{M}\sum_{m=1}^M \frac{\vert \hat\rho_\obs-\hat\rho_m \vert}{\vert \hat\rho_\obs\vert},
\end{align}
where $\Vert \cdot \Vert_1$ denotes the $\ell_1$-norm. 
In order to take care of the randomness involved in sampling the parameters in the matrix $A$, the previous procedure has been replicated 20 times based on 20 independent data samples repetitions. 
The average relative errors $\eta_\K$ and $\eta_\rho$ in~\eqref{eq:relative-errors} were therefore averaged over the 20 independent samples, and reported as $\bar\eta_\K$ and $\bar\eta_\rho$  in Table~\ref{tab:res} (first two rows), along with standard deviations in parentheses.
As for the asymmetry test, we computed for each of the 20 data replications the fraction of times (out of $M$) that the same decision is taken (`reject' vs `do not reject') at the 5\% level, based on the approximate p-values computed on $\bm{X}_m$ and $\bm{X}_\text{obs}$.
 The obtained values were averaged over the 20 independent data samples repetitions and reported in the third row of Table~\ref{tab:res} as $\bar{\text{f}}_{\text{test}}$.

\begin{table}[h!]
\centering
\caption{First two rows: average relative errors~\eqref{eq:relative-errors} for Kendall's distribution function and  Spearman's $\rho$ between the observed sample and the samples retained by the ABC procedure, for varying $K$ (columns).
Third row: fraction of times that the same decision is taken (`reject' vs `do not reject', at the 5\% level) based on $\bm{X}_m$ and  $\bm{X}_\text{obs}$. The results are averaged over 20 independent repetitions. Standard deviations in parentheses. All values are in $\%$.}\label{tab:res}\begin{tabular}{|c||c|c|c|c|}
\hline
{$K$}    & \textbf{2} & \textbf{3} & \textbf{4} & \textbf{5} \\ \hline\hline
$\bar\eta_\K$ & 1.99 (0.16) & 2.32 (0.49) & 2.38 (0.40) & 2.37 (0.56) \\ \hline
$\bar\eta_{\rho}$  & 5.44 (3.98) & 8.16 (5.62) & 6.34 (3.74) & 9.66 (4.12)  \\ \hline
$\bar{\text{f}}_{\text{test}}$ & 16.8 (18.1) & 19.4 (22.0) & 9.2 (13.2) & 13.3 (18.1)  \\  \hline\end{tabular}
\end{table}

Table \ref{tab:res} suggests a general trend: the larger $K$ is, the more difficult the estimation is. However, the estimation procedure yields satisfactory results for all cases considered.

\subsubsection{Misspecified setting: data from a noisy \LF copula\label{sub:Inf_liebFreNoi}}
\begin{figure}[h!]
\centering
 \includegraphics[width=0.3\textwidth]{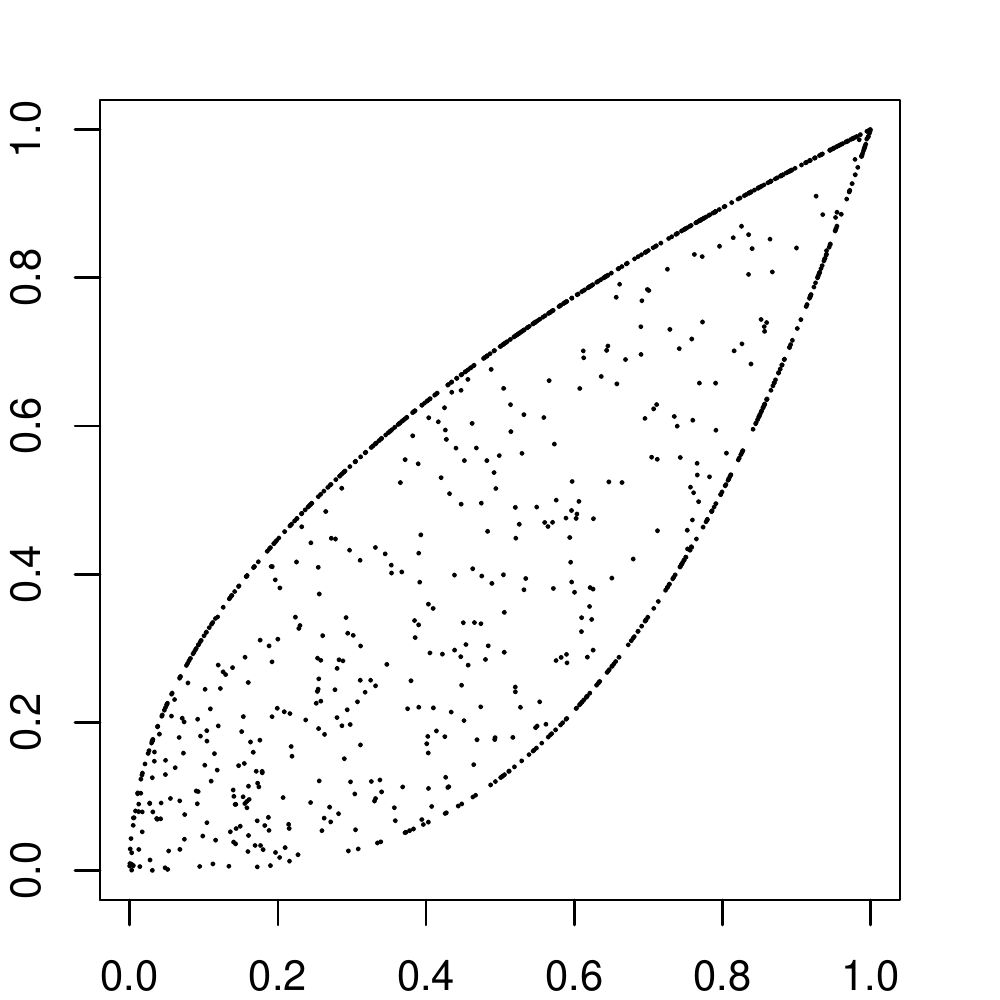}
  \includegraphics[width=0.3\textwidth]{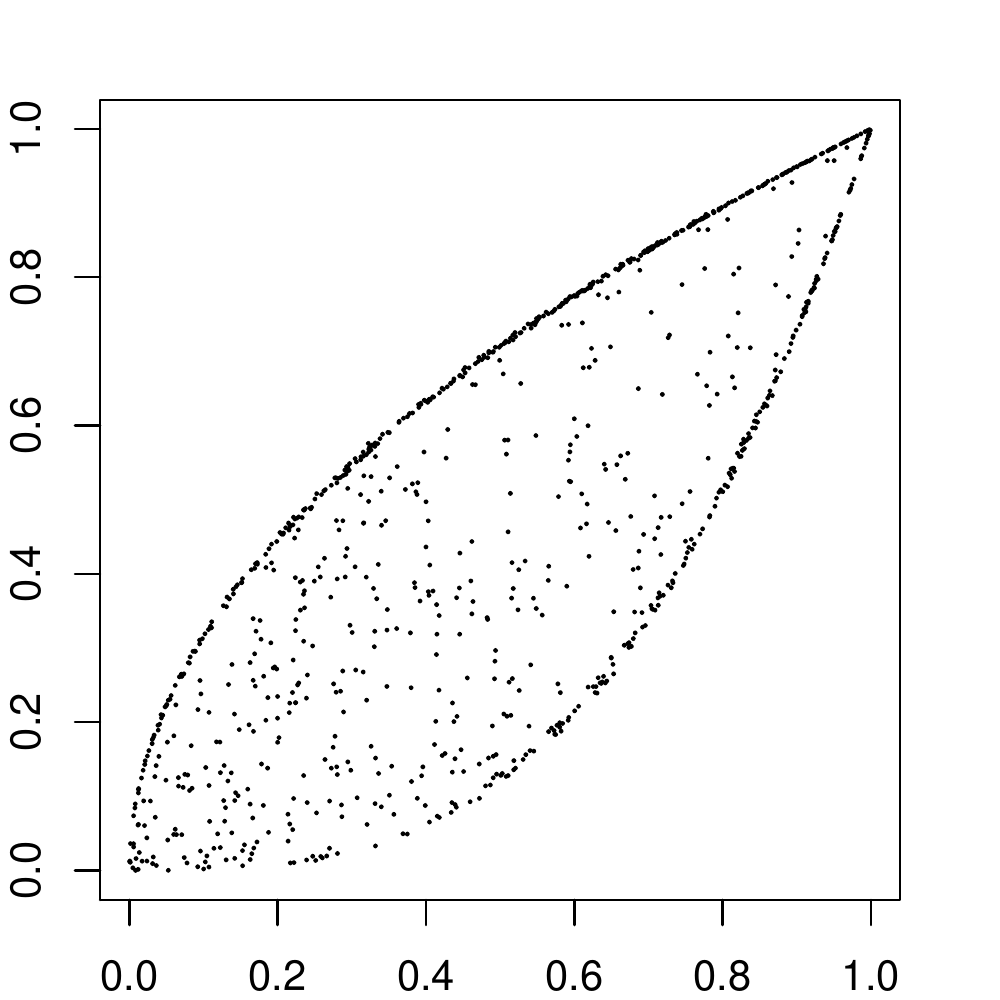}
 \includegraphics[width=0.3\textwidth]{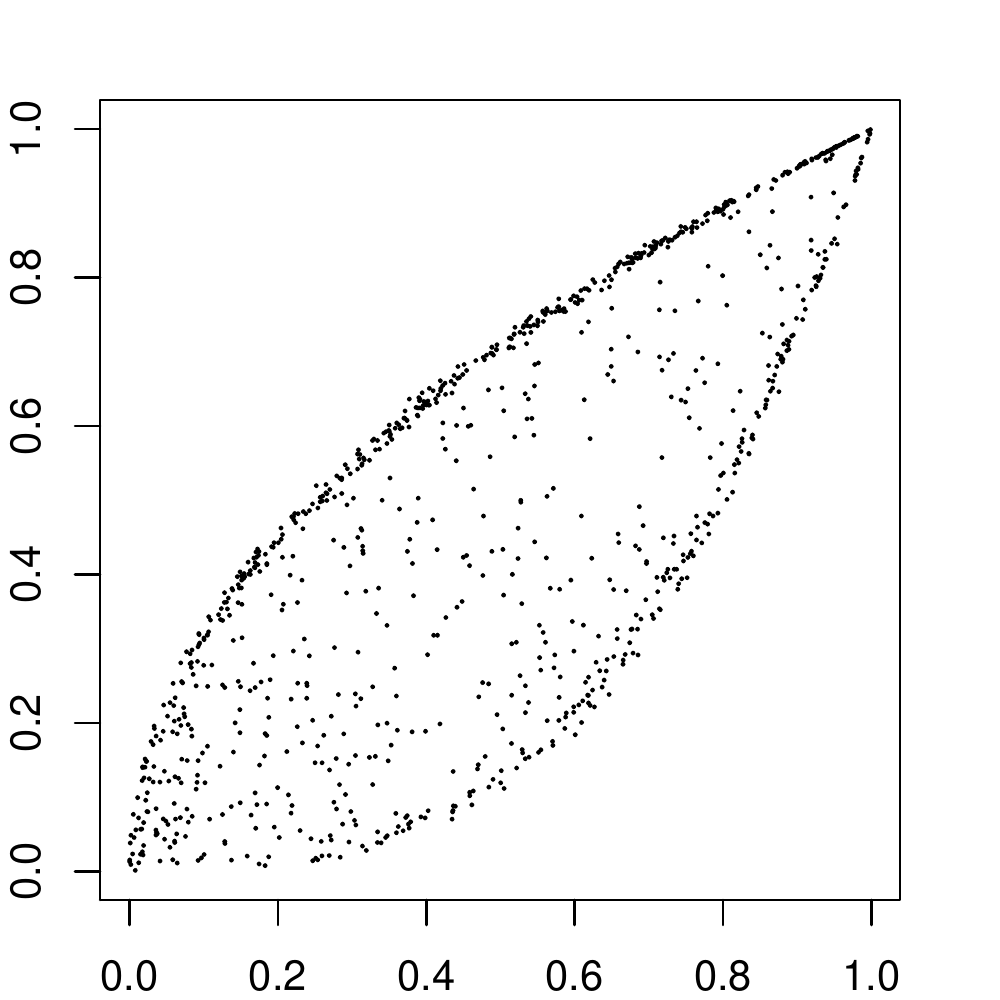}
  \includegraphics[width=0.3\textwidth]{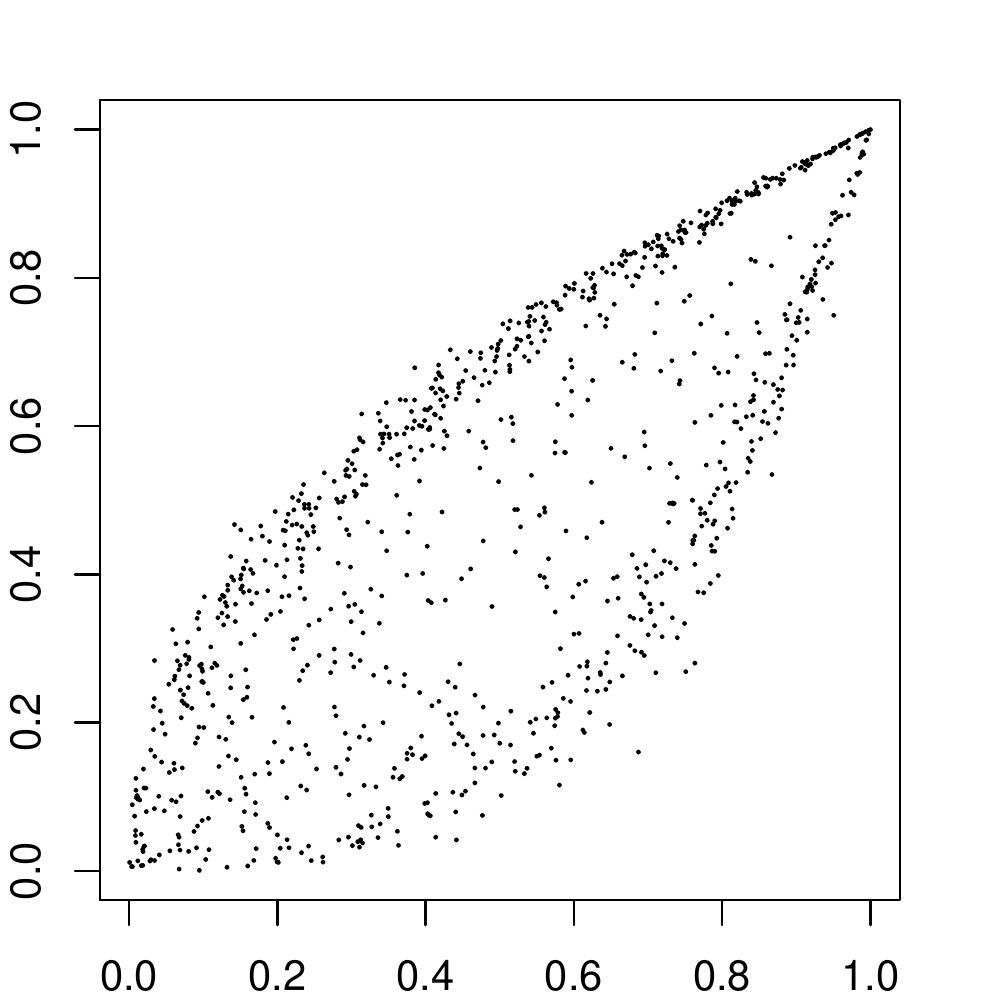}
  \includegraphics[width=0.3\textwidth]{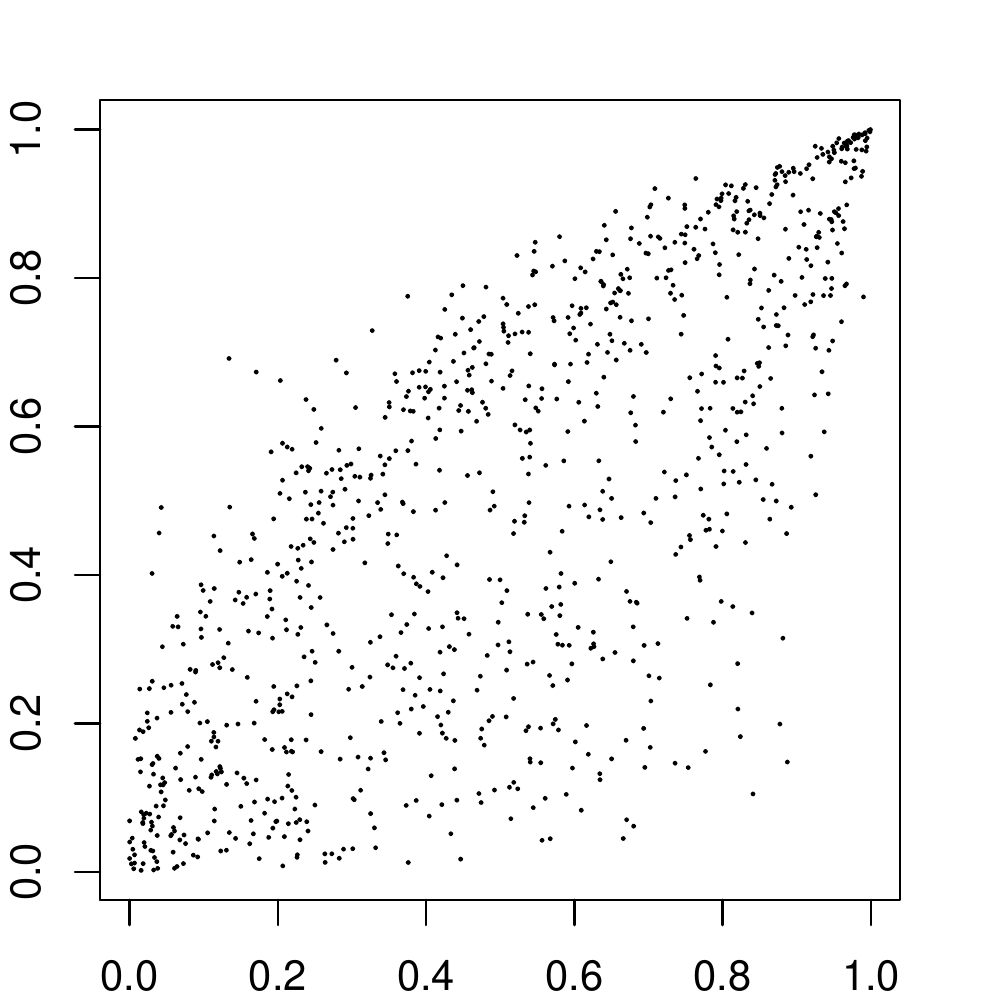}
  \includegraphics[width=0.3\textwidth]{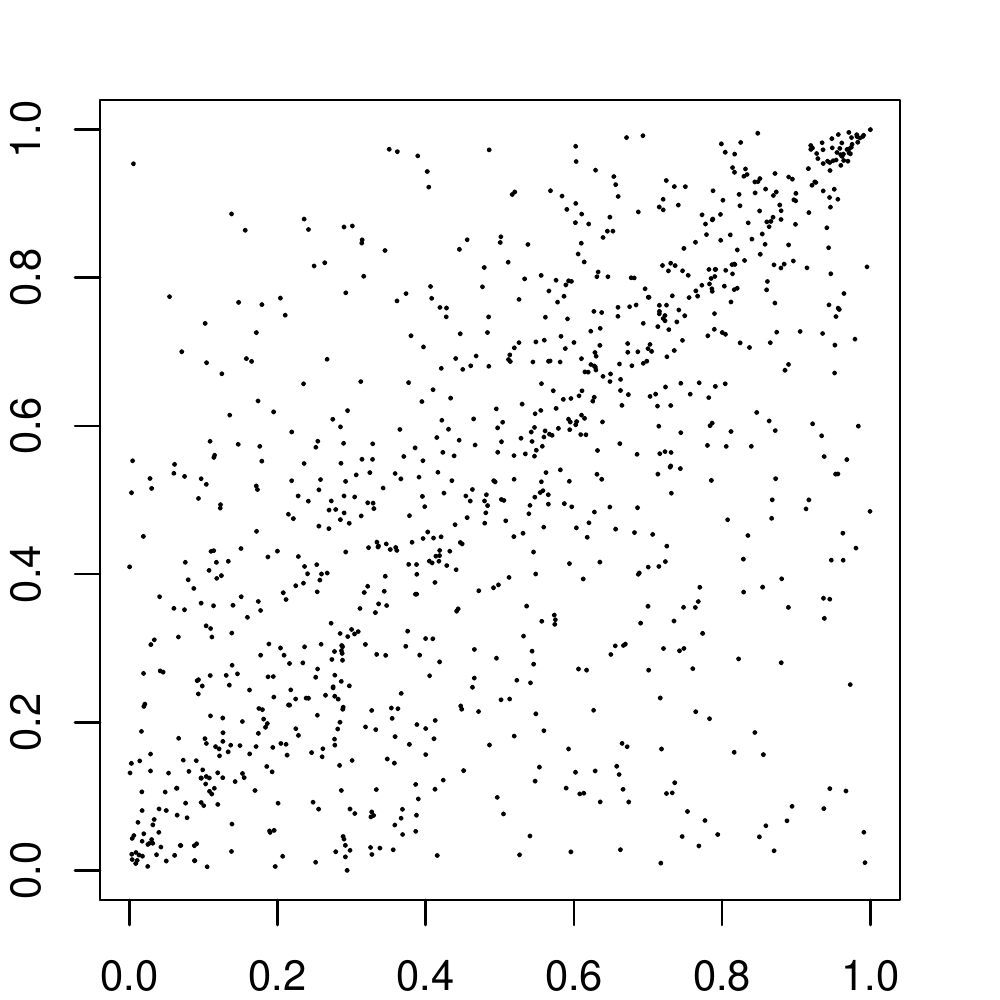}
 \caption{Samples from the exact \LF copula with $d=K=2$, $a_{1}^{(2)}=0.4$, $a_{2}^{(2)}=0.8$ (top-left panel) and from five noisy versions of it, where, for each data point $i\in\{1,\ldots,n\}$, 
 $a_{i1}^{(2)}\sim\mathcal{B}(\alpha_1,\beta_1)$ and $a_{i2}^{(2)}\sim\mathcal{B}(\alpha_2,\beta_2)$  such that ${\rm E}(a_{i1}^{(2)})=0.4$, ${\rm E}(a_{i2}^{(2)})=0.8$, and increasing noise variance $\sigma_a^2 = {\rm var}(a_{i1}^{(2)})={\rm var}(a_{i2}^{(2)})= 10^{-5}, 10^{-4}, 10^{-3}, 10^{-2}, 10^{-1}$  (other panels). Sample size $n=10^3$.}
 \label{fig:noisy}
 \end{figure}

In this section, we generate data from a noisy version of the \LF copula, and demonstrate that our inference procedure works well even if the data are not sampled from the exact model (so-called misspecified setting). In order to sample data from such a noisy model, a slightly changed version of  Algorithm~\ref{algo:sampleFrechet} is used in which the parameters  ${A}$ are not fixed. Instead, they are sampled from a beta distribution with given variance $\sigma_a^2$ (interpreted as the error variance) around some fixed value corresponding to the zero noise version. The latter is illustrated on Fig.~\ref{fig:noisy}: a sample of $n=10^3$ data points from a 2-dimensional \LF copula is depicted on the top-left panel with $K = 2$ iterative steps, and with the two parameters of the power functions set to $a_1^{(2)}=0.4$ and $a_2^{(2)}=0.8$ (recall that $a_1^{(1)}=a_2^{(1)}=1$). The remaining five
panels correspond to samples from \LF copula with increasing noise variance $\sigma_a^2 = 10^{-5}, 10^{-4}, 10^{-3}, 10^{-2}, 10^{-1}$.

The inference results obtained from Algorithm~\ref{algo:abc} (ran with $M'=10^4$ ABC iterations of which $M=300$ were retained) are  reported in Table~\ref{tab:resNoise}.
Unsurprisingly, it appears that the larger the noise variance is, the more difficult is the estimation. Again, the estimation procedure yields good results for all cases considered. 
Let us note that the results reported in Table~\ref{tab:resNoise} are not averaged over 20 independent replications like in the previous section. The reason is that, here, the interest is in illustrating how the procedure deteriorates with the increasing noise in the observed data. 
Therefore, it is sufficient to illustrate the results of the analysis on a single dataset.  

\begin{table}[h!]
\centering
\caption{First two rows: average relative errors~\eqref{eq:relative-errors} for Kendall's distribution function and  Spearman's $\rho$ between the observed sample and the samples retained by the ABC procedure, for growing noise (columns).
Third row: fraction of times that the same decision is taken at the 5\% level based on $\bm{X}_m$ and $\bm{X}_\text{obs}$. All values are in $\%$.}\label{tab:resNoise}
\begin{tabular}{|c||c|c|c|c|c|c|}
\hline
$\sigma_a^2$    & 0 & $10^{-5}$ & $10^{-4}$ & $10^{-3}$ & $10^{-2}$ & $10^{-1}$ \\ \hline\hline
$\eta_\K$ & 1.7 & 1.68 & 1.77 & 1.68 & 2.03 & 2.81\\ \hline
$\eta_\rho$  & 4.07 & 3.77 & 4.66 & 4.30 & 4.88 & 13.00 \\ \hline
$\text{f}_{\text{test}}$ & 9.00 & 1.33 & 15.7 & 3.33 & 2.33 & 7.00 \\\hline
\end{tabular}
\end{table}



{\subsection{Comparison of ABC with likelihood-based estimation\label{sub:abc_VS_mle}}
In this section, the inference based on our ABC procedure is compared with the likelihood-based method provided in the \texttt{copula}  \textsf{R} package \citep{copulaR}. 

As a matter of fact, in some specific cases it is possible to derive the density of Liebscher copulas, and hence, to perform likelihood-based inference on them. We work here with the following bivariate 
Liebscher copula obtained by combining Clayton copula $C_1(u,v)=\left(u^{-\theta}+v^{-\theta}-1\right)^{-1/\theta}$ with the independence copula  $C_2(u,v)=uv$, and by using power functions (Example~\ref{expower}),
\begin{equation}\label{compare}
C(u,v)=C_1(u^{p}, v^{q})C_2(u^{1 - p}, v^{1 - q}).
\end{equation}
In this illustration we set $\theta = 5$, $
p = 0.3$, and $q = 0.8$, and estimate these parameters on 100 data sets independently sampled from copula~\eqref{compare}, with both our ABC procedure and the optimization procedure implemented in the \texttt{fitCopula}  function of the \texttt{copula} package. 
The above procedure is repeated for samples of size $n=500$ and $n=10\,000$. In order to compare the results of our Bayesian procedure with the likelihood-based one, which provides only point estimates (maximum likelihood estimator, MLE), the posterior distribution of the model parameters is summarized into two point estimates: the posterior mean and the posterior median. 

The results are plotted in Fig.~\ref{fig:compare}, where each boxplot is made out of the 100 estimated values. In each panel, three quantities are plotted: \texttt{MLE}, which refers to the likelihood-based estimation, \texttt{Post.median}, which refers to the posterior median ABC estimation, and \texttt{Post.mean}, which refers to the posterior mean ABC estimation. It appears from the three plots that, as expected, the MLE is asymptotically unbiased. This can be seen by noticing the convergence to the true value as $n$ increases.
On the other hand, both the mean and the median of the posterior distribution obtained with our ABC procedure show some bias. Unexpectedly, the results from MLE show a larger variance than ABC in some cases, particularly when $n=500$. This could be due to difficulties in the optimisation procedure of the \textsf{fitCopula} \textsf{R} function which, for small sample sizes, could have problems converging.

\begin{figure}[h!]
\centering
        \subfloat[Parameter $p$ (first coordinate)]
        {\includegraphics[width=.3\textwidth]{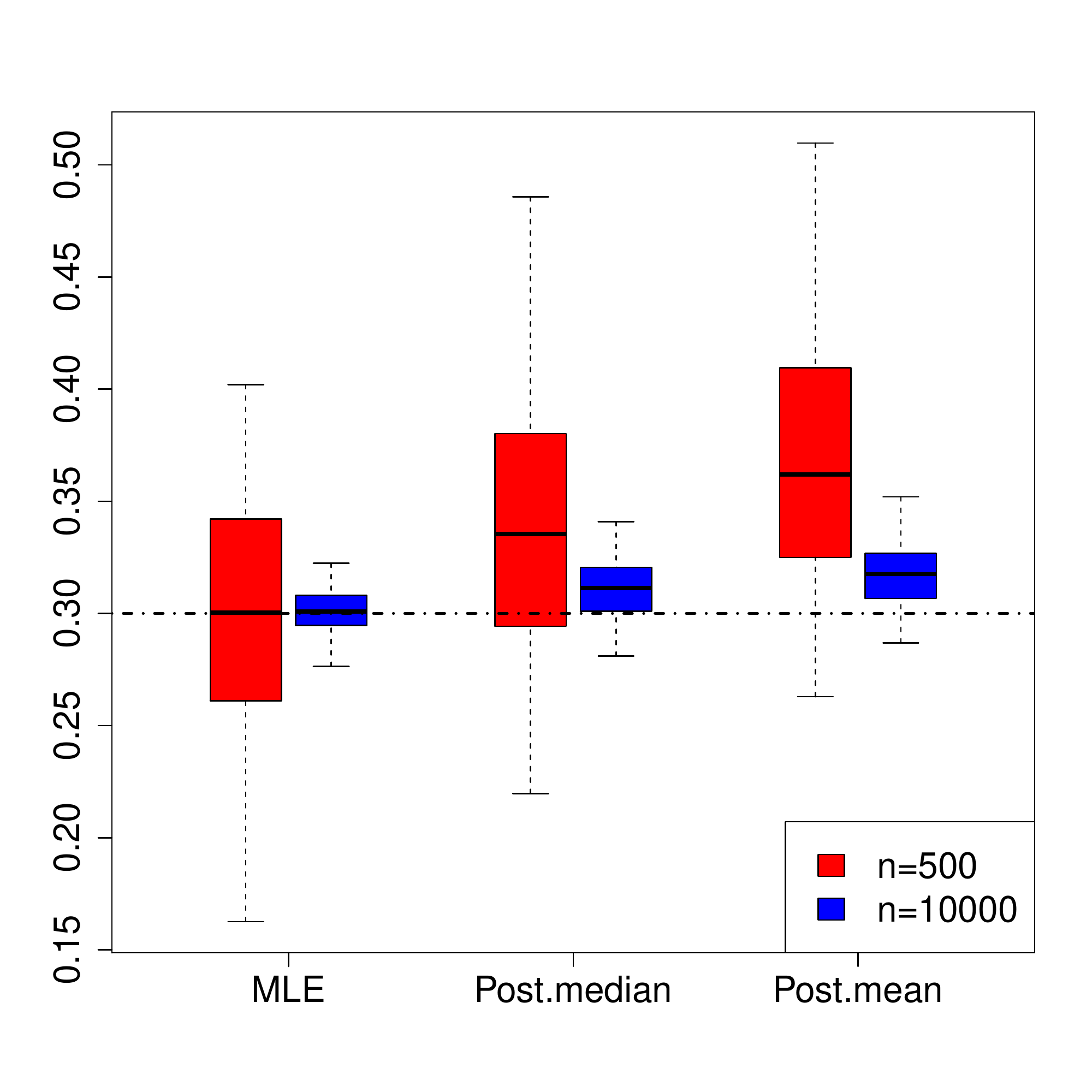}}
        ~
        \subfloat[Parameter $q$ (second coordinate)]
        {\includegraphics[width=.3\textwidth]{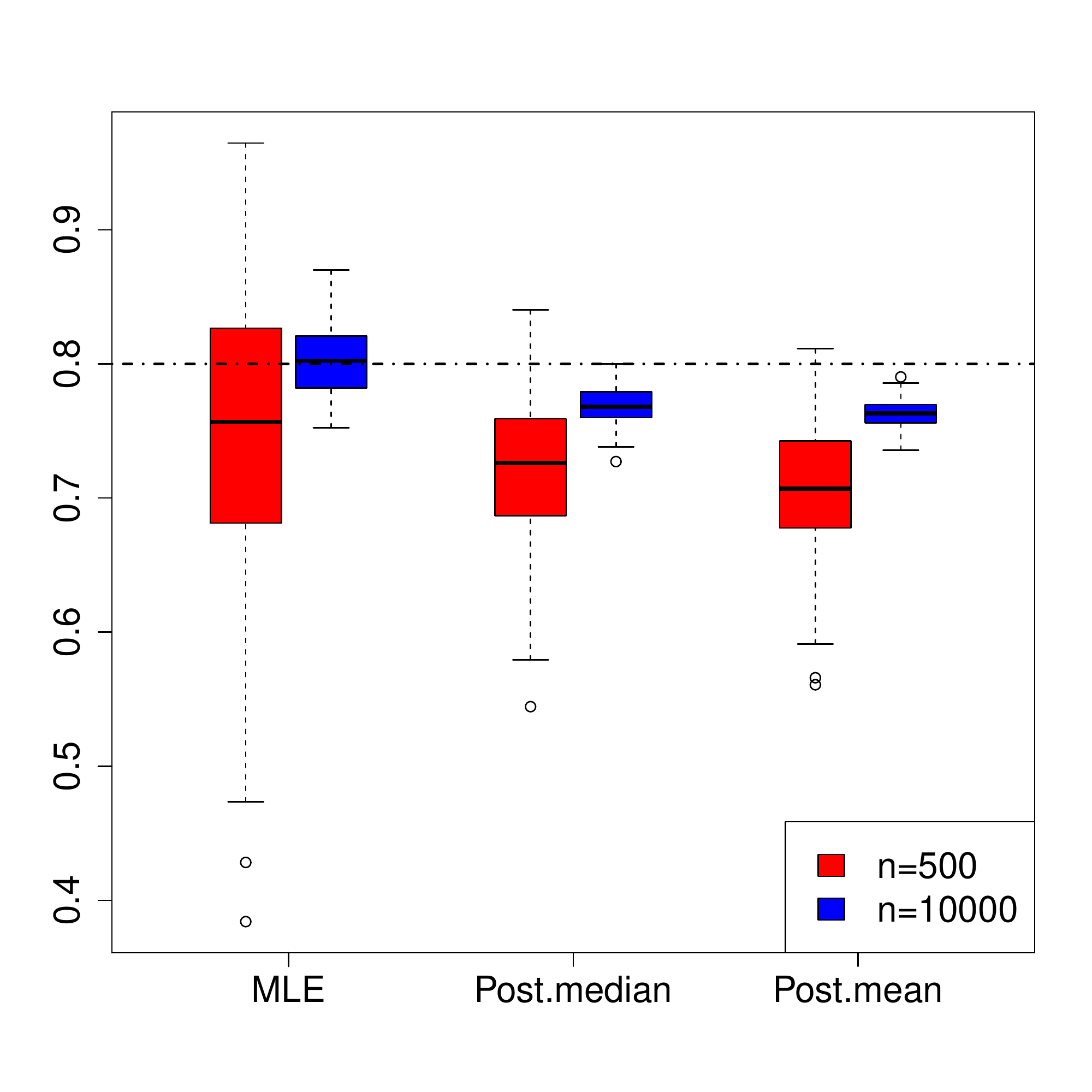}}
        ~
        \subfloat[Parameter $\theta$ (Clayton copula)]
        {\includegraphics[width=.3\textwidth]{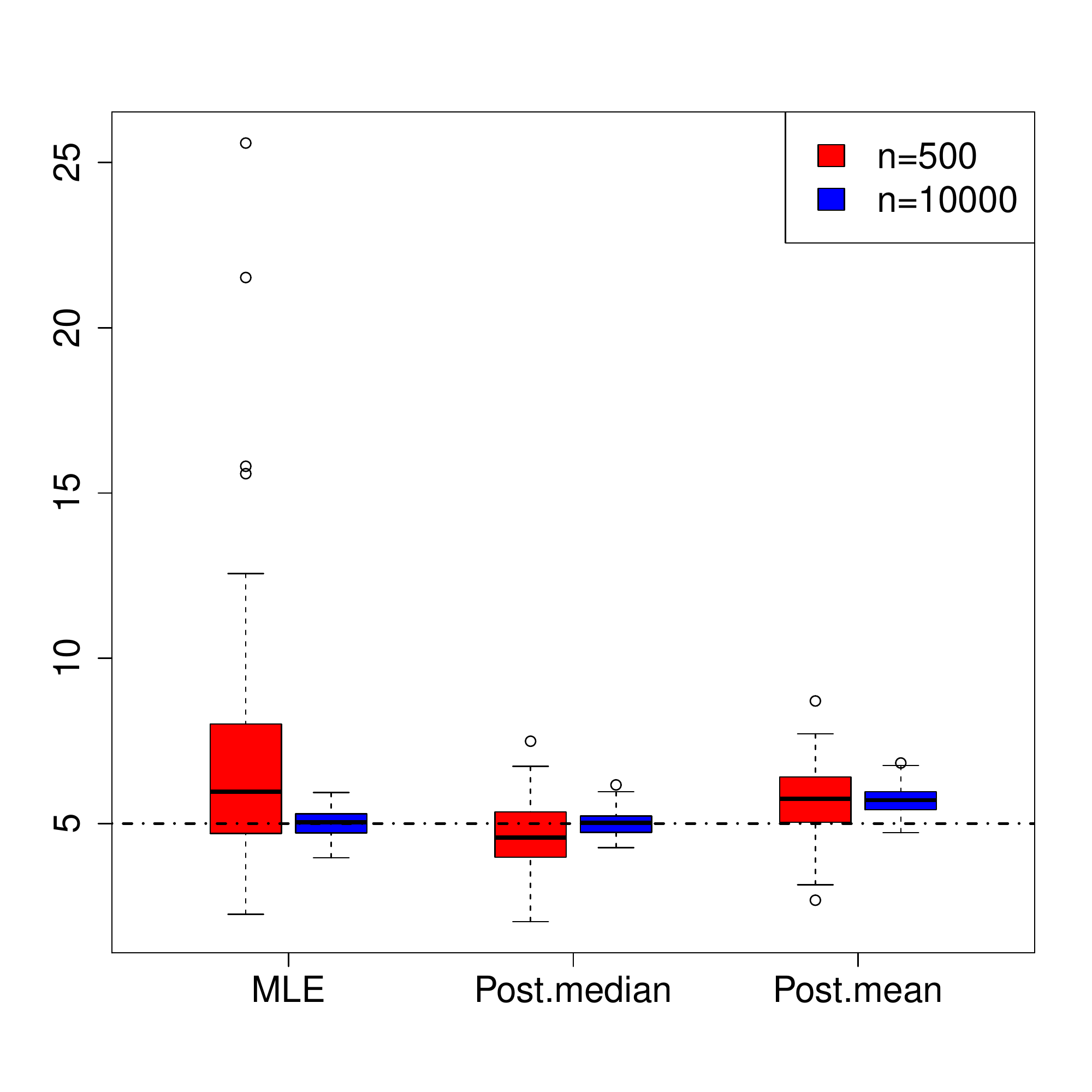}}
\caption{Results from the comparison. (a): boxplot of the parameter $p$; (b): boxplot of the parameter $q$; (c): boxplot of the parameter $\theta$. Red  corresponds to samples of size $n=500$, blue to samples of size $n=10\,000$. The horizontal dotted line represents  the true value.} \label{fig:compare}
\end{figure}

}

\section{Conclusion}\label{sec:concl}
In this paper, we have studied the class of asymmetric copulas first introduced by \citet{Khoudraji}, and developed in its general form by \citet{Liebscher}. Some new theoretical properties of these copulas were provived, including novel closed form expressions for its tail dependence coefficients, thus complementing the partial results of \citet{Liebscher} and \citet{Liebscher2}. An iterative procedure is also introduced to flexibly sample from these copulas, which makes it easy to apply an Approximate Bayesian computation procedure to make inference on them. 


\section*{Acknowledgments}
The authors would like to thank two referees and the Associate Editor for their valuable suggestions,
which have significantly improved the paper. 
The work of the third author benefited from the support of the Chair Stress Test, RISK Management and Financial Steering, led by the French \'Ecole polytechnique and its Foundation and sponsored by BNP Paribas.

\appendix

\section{Proofs of main results}

\paragraph{\bf Proof of Proposition~\ref{prop:L}}

(i) Recall that $d=2$ and let  $g^{(k)}\coloneqq g_1^{(k)}=g_2^{(k)}$ for all $k\in\{1,\dots,K\}$. For all $x>0$, we have
\begin{align*}
x = \frac{\varepsilon x}{\varepsilon} \overset{\text{(a)}}{=} \prod_{k=1}^K \frac{g^{(k)}(\varepsilon x)}{g^{(k)}(\varepsilon)} 
\stackrel[\varepsilon\to 0]{\text{(b)}}{\longrightarrow}
\prod_{k=1}^K x^{\gamma^{(k)}} = x^{\sum_{k=1}^K \gamma^{(k)}},
\end{align*}
since the product of the $g^{(k)}$ functions is the identity (a) and by definition of regular variation (b). It follows that $\sum_{k=1}^K \gamma^{(k)}=1$. 
Besides,
\begin{align*}
\frac{\tilde{C}(\varepsilon x, \varepsilon y)}{\varepsilon} 
&= \prod_{k=1}^K \frac{C_k(g^{(k)}(\varepsilon x),g^{(k)}(\varepsilon y))}{g^{(k)}(\varepsilon)} \overset{\text{(a)}}{=} 
\prod_{k=1}^K \frac{C_k(g^{(k)}(\varepsilon) x^{\gamma^{(k)}}+o(g^{(k)}(\varepsilon)),g^{(k)}(\varepsilon)y^{\gamma^{(k)}}+o(g^{(k)}(\varepsilon)))}{g^{(k)}(\varepsilon)} \\
&\overset{\text{(b)}}{=}
\prod_{k=1}^K \frac{C_k(g^{(k)}(\varepsilon) x^{\gamma^{(k)}},g^{(k)}(\varepsilon)y^{\gamma^{(k)}})+o(g^{(k)}(\epsilon))}{g^{(k)}(\varepsilon)} 
=
\prod_{k=1}^K \frac{C_k(g^{(k)}(\varepsilon) x^{\gamma^{(k)}},g^{(k)}(\varepsilon)y^{\gamma^{(k)}})}{g^{(k)}(\varepsilon)} +o(1)\\
&\stackrel[\varepsilon\to 0]{\text{(c)}}{\longrightarrow}
\prod_{k=1}^K \Lambda_L(C_k; x^{\gamma^{(k)}}, y^{\gamma^{(k)}}),
\end{align*}
by the regular variation property (a), by Lipschitz property for copulas and because of the continuity of
$g^{(k)}$ at the origin (c). The result is thus proved. \\
(ii) Let $x\leq y$, the proof being similar when $x>y$. We have
\begin{align*}
\tilde{C}(\varepsilon x,\varepsilon y) = \prod_{k=1}^K C_k(g_1^{(k)}(\varepsilon x),g_2^{(k)}(\varepsilon y) ){\leq} 
\prod_{k=1}^K \min(g_1^{(k)}(\varepsilon x),g_2^{(k)}(\varepsilon y)) ,
\end{align*}
since, for any copula $C$ and any $(u,v)\in[0,1]^2$, $C(u,v)\leq\min(u,v)$. By assumption here, there exists
$k_0\in\{1,\dots,K\}$ such that $\eta(t)\coloneqq g_1^{(k_0)}(t)/g_2^{(k_0)}(t)$, where $\eta(t)\to 0$ as $t\to 0$. 
Taking into account that $g_1^{(k_0)}$ is increasing, this entails
\begin{align*}
\prod_{k=1}^K \min(g_1^{(k)}(\varepsilon x),g_2^{(k)}(\varepsilon y)) \leq g_1^{(k_0)}(\varepsilon x) \prod_{k\neq k_0}g_2^{(k)}(\varepsilon y) \leq g_1^{(k_0)}(\varepsilon y) \prod_{k\neq k_0}g_2^{(k)}(\varepsilon y) = \eta(\varepsilon y) \prod_{k=1}^K g_2^{(k)}(\varepsilon y).
\end{align*}
Recalling that $\prod_{k=1}^K g_2^{(k)}$ is the identity function, it follows that
$
\tilde{C}(\varepsilon x,\varepsilon y) \leq \eta(\varepsilon y) \varepsilon y
$
and consequently $\Lambda_L(\tilde{C}; x,y) = 0$.\\
(iii) Differentiating the product of the $g_j^{(k)}$ functions, we obtain for $j\in\{1,2\}$:
\begin{align*}
1=\sum_{k=1}^K\big(g_j^{(k)}\big)^\prime (u) \prod_{\ell\neq k} g_j^{(\ell)}(u).
\end{align*}
Now, since $g_j^{(\ell)}(1)=1$ for all $\ell\in\{1,\ldots,K\}$, and $j=1,2$, we obtain  $\sum_{k=1}^K \text{d}_j^{(k)}=1$. 
Turning to the upper dependence function, let us write
\begin{align}\label{eq:upperL_1}
\tilde{C}(1-\varepsilon x, 1-\varepsilon y)-1 \stackrel[\varepsilon\to 0]{}{\sim} - \ln \tilde{C}(1-\varepsilon x, 1-\varepsilon y)  = -\sum_{k=1}^K\ln C_k(g_1^{(k)}(1-\varepsilon x), g_2^{(k)}(1-\varepsilon y)),
\end{align}
where $ \stackrel[\varepsilon\to 0]{}{\sim}$ denotes the asymptotic equivalence as $\varepsilon\to 0$.
By assumption, all $g_j^{(k)}$ are differentiable at 1, with derivative denoted by $\text{d}_j^{(k)}$, and satisfy $g_j^{(k)}(1)=1$. Hence, the first order Taylor expansion is $g_j^{(k)}(1-\varepsilon) = 1- \text{d}_j^{(k)}\varepsilon+o(\varepsilon)$. By the Lipschitz property of copulas,
$$
C_k(g_1^{(k)}(1-\varepsilon x), g_2^{(k)}(1-\varepsilon y)) 
= C_k(1- \text{d}_1^{(k)}\varepsilon x+o(\varepsilon), 1- \text{d}_2^{(k)}\varepsilon y+o(\varepsilon)) 
=
C_k(1- \text{d}_1^{(k)}\varepsilon x, 1- \text{d}_2^{(k)}\varepsilon y)+o(\varepsilon). 
$$
Plugging in to~\eqref{eq:upperL_1} yields
\begin{align*}
\frac{\tilde{C}(1-\varepsilon x, 1-\varepsilon y) -1}{\varepsilon}
&\stackrel[\varepsilon\to 0]{}{\sim} 
-\frac{1}{\varepsilon}\sum_{k=1}^K\ln\left(C_k(1- \text{d}_1^{(k)}\varepsilon x, 1- \text{d}_2^{(k)}\varepsilon y)+o(\varepsilon)\right)\\
&\stackrel[\epsilon\to 0]{}{\sim} 
\sum_{k=1}^K \frac{C_k(1- \text{d}_1^{(k)}\varepsilon x, 1- \text{d}_2^{(k)}\varepsilon y)-1}{\varepsilon}+o(1).
\end{align*}
Finally, in view of $\sum_{k=1}^K \text{d}_j^{(k)}=1$ for $j=1,2$ it follows that
\begin{align*}
x+y+\frac{\tilde{C}(1-\varepsilon x, 1-\varepsilon y) -1}{\varepsilon}
&\stackrel[\varepsilon\to 0]{}{\sim} 
\sum_{k=1}^K \left(d_1^{(k)} x + d_2^{(k)} y+\frac{C_k(1- \text{d}_1^{(k)}\varepsilon x, 1- \text{d}_2^{(k)}\varepsilon y)-1}{\varepsilon}\right)+o(1).
\end{align*}
Taking the limit $\varepsilon\to 0$ yields the result: 
$$
\Lambda_U(\tilde{C}; x,y) = \sum_{k=1}^K \Lambda_U(C_k; \text{d}_1^{(k)}x,\text{d}_2^{(k)} y).
$$
(iv) The case $\text{d}_1^{(k)}=\text{d}_2^{(k)}$ for all $k\in\{1,\dots,K\}$ is then a simple consequence of the homogeneity property of the upper tail dependence function.
\CQFD


\paragraph{\bf Proof of Corollary~\ref{cor:lambda}}  This is a direct consequence of Proposition~\ref{prop:L}. \CQFD


\paragraph{\bf Proof of Proposition~\ref{prop:long_list_of_dependence}}  

\citet[Proposition 2.2]{Liebscher} proves the \TPtwo and  LTD properties. The LTI property can be proven in the same way as LTD. 
For PQD, it suffices to remark that for any $u,v\in[0,1]$, 
\begin{align*}
\tilde{C}(u,v) = \prod_{k=1}^K C_k(g_1^{(k)}(u),g_2^{(k)}(v)) 
\overset{\text{(a)}}{\geq}  \prod_{k=1}^K g_1^{(k)}(u)g_2^{(k)}(v)
= \prod_{k=1}^K g_1^{(k)}(u)  \prod_{k=1}^K g_2^{(k)}(v) = uv,
\end{align*}
while NQD works similarly with a reversed inequality in $(a)$.

Let us prove the SI part. According to Equation~\eqref{eq:SI-charact}, \SI is a property of $u\mapsto \tilde{C}(u,v)$  and $v\mapsto \tilde{C}(u,v)$  functions. Focusing without loss of generality on the former function, and omitting the $v$ variable for notational simplicity,  $\tilde{C}$ can be written as
\begin{align*}
\tilde{C} (u) = C_1(g^{(1)}(u))\ldots C_K(g^{(K)}(u)),
\end{align*}
with $g^{(1)}(u)\ldots g^{(K)}(u)= u$ for all $u\in[0,1]$.
Differentiating this function twice yields $\tilde{C}''(u)=(T_1(u)+T_2(u))\tilde{C}(u)$ where
\begin{align*}
T_1 =  \sum_{k=1}^K (g^{(k)\prime})^2
\frac{C_k''(g^{(k)})}{C_k(g^{(k)})}, \quad
T_2 = \sum_{k=1}^K \tau_k \frac{C_k'(g^{(k)})}{C_k(g^{(k)})}, 
\quad 
\tau_k = g^{(k)\prime\prime}+g^{(k)\prime}
\sum_{\ell\neq k}g^{(\ell)\prime}\frac{C_\ell'(g^{(\ell)})}{C_\ell(g^{(\ell)})}.
\end{align*}
By assumption, $C_1,\ldots,C_K$ are SI, hence they are concave, thus $C_1''\leq 0$,\ldots, $C_K''\leq 0$ and therefore $T_1\leq 0$.  In view of Theorem~5.2.12 and Corollary~5.2.6 in~\citet{Nelsen}, $u\mapsto C_\ell'(u)-{C_\ell(u)}/{u}$ is a negative function for all $\ell\in\{1,\ldots,K\}$. As a consequence, $\tau_k$ can be upper bounded as follows:
\begin{align}\label{eq:tk-bound}
    \tau_k\leq g^{(k)\prime\prime}+g^{(k)\prime}
\sum_{\ell\neq k}\frac{g^{(\ell)\prime}}{g^{(\ell)}}
{=} g^{(k)\prime\prime}+g^{(k)\prime}
\left(\frac{1}{\Id} -\frac{g^{(k)\prime}}{g^{(k)}}\right),
\end{align}
where the equality is due to the fact that the product of all $g^{(\ell)}$ functions is the identity $\Id$, thus the derivative of product logarithm is $1/\Id$. Additionally, since $g^{(k)}$ is concave and $g^{(k)}(0)=0$, Theorem 5 by \cite{bruckner1962some} implies that $-g^{(k)}$ is star-shaped, i.e., $-g^{(k)}/Id$
is increasing. This, in turn, proves that the right-hand side of~\eqref{eq:tk-bound} can be further upper bounded by zero, 
and therefore $T_2\leq 0$.
As a conclusion, $\tilde{C}''<0$ and, in virtue of Equation~\eqref{eq:SI-charact}, $\tilde{C}$ is also SI. The proof for SD follows similar lines.
\CQFD


\paragraph{\bf Proof of Proposition~\ref{propstable1}} 

By definition of the copula $\tilde{C}$,
\begin{align*}
 \ln \tilde{C}(\bm{u})&=  \sum_{k=1}^K \sum_{i=1}^m \theta_{ik} \ln\varphi_i\left(u_1^{p^{(k)}},\dots,u_d^{p^{(k)}}\right)
=\sum_{k=1}^K \sum_{i=1}^m \theta_{ik} 
(\ln\mathop \circ\varphi_i\mathop \circ\exp) \left(p^{(k)}\ln u_1,\dots,p^{(k)} \ln u_d\right)\\
&=\sum_{k=1}^K \sum_{i=1}^m \theta_{ik} (p^{(k)})^{\lambda_i}
(\ln\mathop \circ\varphi_i\mathop \circ\exp) \left(\ln u_1,\dots, \ln u_d\right)
=\sum_{k=1}^K \sum_{i=1}^m \theta_{ik} (p^{(k})^{\lambda_i}
\ln\varphi_i \left( \bm{u}\right)
=\sum_{i=1}^m \sum_{k=1}^K  \theta_{ik} (p^{(k)})^{\lambda_i}
\ln\varphi_i \left( \bm{u}\right) \\
&=\sum_{i=1}^m   \tilde\theta_{iK} \ln\varphi_i \left( \bm{u}\right) 
= \ln  C(\bm{u} \;|\; \tilde \theta_{1K},\dots,\tilde\theta_{mK}),
\end{align*}
and the result is proved.\CQFD


\paragraph{\bf Proof of Proposition~\ref{propstable2}}

Let us remark that Example~\ref{exMaxStable} shows that
$C$ is max-stable implies $\tilde{C}^{(K)}=C$ for all $K\geq 1$ and for all sequence  $(p^{(k)})_k\subset(0,1)$.
Conversely, assume that $\tilde{C}^{(K)}=C$ for all $K\geq 1$ and for all sequence  $(p^{(k)})_k\subset(0,1)$.
From~(\ref{eq:lieb_cop}) with $K=2$, it follows that
$$
C(\bm{u})=C\left(u_1^{p^{(1)}},\dots,u_d^{p^{(1)}}\right)C\left(u_1^{1-p^{(1)}},\dots,u_d^{1-p^{(1)}}\right),
$$
for all $\bm{u}\in[0,1]^d$. Introducing $\varphi:{\mathbb R}_-^d\to {\mathbb R}_-$ the continuous function 
defined by $\varphi=\ln\mathop \circ C\mathop \circ \exp$, we thus have
 $\varphi(\bm{v}) = \varphi(p^{(1)} \bm{v}) + \varphi((1-p^{(1)}) \bm{v})$,
for all $\bm{v}\in{\mathbb R}_-^d$ and $p^{(1)}\in(0,1)$. Lemma~\ref{lemhomo} in~\ref{appendix:aux} entails that 
$\varphi$ is homogeneous of degree 1 or equivalently that $C$ is max-stable.\CQFD


\paragraph{\bf Proof of Proposition \ref{prop-rewrite}}

Let us first show that the copula $\tilde{C}^{(K)}$, $K\geq 1$ defined iteratively by~(\ref{def1}),~(\ref{def2}) is a Liebscher copula.
The proof is done by induction on $K$.
First, it is clear that
$$
\tilde{C}^{(1)}(\bm{u}) = C_1\left(g_1^{(1,1)}(u_1),\dots,g_d^{(1,1)}(u_d) \right) = 
C_1\left(u_1/f^{(1)}(u_1),\dots,u_d/f^{(1)}(u_d) \right)
= C_1(\bm{u}).
$$
Secondly, assume that
$$
  \tilde{C}^{(K-1)}(\bm{u}) = \prod_{k=1}^{K-1} 
C_k\left(g_1^{(K-k,K-1)}(u_1),\dots,g_d^{(K-k,K-1)}(u_d) \right)
$$
and
\begin{align}
\label{hyp1}
 g_j^{(1,K-1)}&=\Id/f_j^{(K-1)}, \\
 \label{hyp2}
g_j^{(k,K-1)}&= \bigodot_{i=K-k-1}^{K-1} f_j^{(i)} \left/ \bigodot_{i=K-k}^{K-1} f_j^{(i)} \right. \; , k\in\{2,\dots,K-1\},
\end{align}
where the $\bigodot$ notation is defined in~(\ref{notation-rond}).
From (\ref{def2}), it follows that
\begin{align}
\nonumber
 \tilde{C}^{(K)}(\bm{u})&= C_K\left(\frac{u_1}{f^{(K)}_1(u_1)},\dots,\frac{u_d}{f^{(K)}_d(u_d)} \right)
 \tilde{C}^{(K-1)}\left(f^{(K)}_1(u_1),\dots,f^{(K)}_d(u_d)\right) \\
 \nonumber
&= C_K\left(\frac{u_1}{f^{(K)}_1(u_1)},\dots,\frac{u_d}{f^{(K)}_d(u_d)} \right) 
 \prod_{k =1}^{K-1} C_k \left(g_1^{(K-k,K-1)}\mathop \circ f^{(K)}_1(u_1),\dots,g_d^{(K-k,K-1)}\mathop \circ 
f^{(K)}_d(u_d)\right) \\
\label{eqtmp0}
 &= \prod_{k=1}^K C_k\left(g_1^{(K-k+1,K)}(u_1),\dots,g_d^{(K-k+1,K)}(u_d) \right)
\end{align}
by letting for all $j\in\{1,\dots,d\}$:
\begin{align}
\label{eqtmp1}
g_j^{(K-k+1,K)}&= g_j^{(K-k,K-1)}\mathop \circ f^{(K)}_j, \; 1 \leq k \leq K-1, \\
\label{eqtmp2}
g_j^{(1,K)} &= \Id / f_j^{(K)}.
\end{align}
As a first result, (\ref{eqtmp0}) proves~(\ref{def3}) while~(\ref{eqtmp2}) proves~(\ref{gj1}).
Letting $\ell=K-k+1$ and in view of~(\ref{hyp2}), equation~(\ref{eqtmp1}) can be rewritten
 for all  $2 \leq \ell \leq K$ as
$$
g_j^{(\ell,K)}= g_j^{(\ell-1,K-1)}\mathop \circ f^{(K)}_j 
 = \left(  \bigodot_{i=K-\ell+2}^{K-1} f_j^{(i)} \left/ \bigodot_{i=K-\ell+1}^{K-1} f_j^{(i)} \right.\right) \mathop \circ  f^{(K)}_j 
=  \bigodot_{i=K-\ell+2}^{K} f^{(i)}_j \left/  \bigodot_{i=K-\ell+1}^{K} f^{(i)}_j \right.
$$
which proves~(\ref{gj2}) for $\ell\in\{2,\dots,K\}$. The case $\ell=1$ is straightforward in view of~(\ref{eqtmp2}).\\
Conversely, let us prove that any Liebscher copula can be constructed iteratively by~(\ref{def1}),~(\ref{def2}).
Let $g_j^{(k)}$, $j\in\{1,\dots,d\}$, $k\in\{1,\dots,K\}$ be the set of functions associated with Liebscher  copula~(\ref{eq:lieb_cop}).
Our goal is to find a set of functions $f_j^{(k)}$, $j\in\{1,\dots,d\}$, $k\in\{1,\dots,K\}$ verifying the set of equations~(\ref{def1}),~(\ref{def2}).
For all $j\in\{1,\dots,d\}$, define $f_j^{(K)}\coloneqq \Id/g_j^{(1)}$ 
and, for 
all $k\in\{K-1,\dots,1\}$, 
\begin{equation}
 \label{Gj1}
f_j^{(k)}\coloneqq  \Id\left/ \left( g_j^{(K+1-k)} \mathop \circ(F_j^{(k+1)})^{-1} \right.\right) \mbox{ where } 
F_j^{(k+1)}\coloneqq   \bigodot_{i=k+1}^{K} f_j^{(i)}.
\end{equation}
Let $j\in\{1,\dots,k\}$. The first part of the proof consists in establishing that
$F_j^{(k+1)}$ is strictly increasing $[0,1]\to[0,1]$
and thus that its inverse $(F_j^{(k+1)})^{-1}$ is well-defined.
To this end, remark that $F_j^{(K)}=\Id$ by definition and $F_j^{(k)}=  f_j^{(k)} \mathop \circ F_j^{(k+1)}= F_j^{(k+1)}/g_j^{(K+1-k)}$,
in view of~(\ref{gj2}).
Iterating, it follows that 
$$
F_j^{(k+1)}=\Id \left/ \prod_{i=2}^{K-k} g_j^{(i)} \right.= g_j^{(1)} \prod_{i=K-k+1}^{K} g_j^{(i)}
$$
in view of (\ref{eq:lieb_assumption}). It is then clear that $F_j^{(k+1)}$ is strictly increasing $[0,1]\to[0,1]$.
The goal of the second part of the proof is to show that 
$f_j^{(k)} \in{\cal F}$ for
all $k\in\{1,\dots,K\}$. It is clear 
that $\Id/f_j^{(k)}$ is strictly increasing for all $k\in\{1,\dots,K\}$ from~(\ref{Gj1}).
Besides, as already noticed, $f_j^{(k)} \mathop \circ F_j^{(k+1)} = F_j^{(k)}$,
and thus $f_j^{(k)}$
is strictly increasing as the composition of strictly increasing functions.  This concludes the proof that $f_j^{(k)} \in{\cal F}$.
The third part of the proof consists in showing that $f_j^{(k)}$, $k\in\{1,\dots,K\}$ is solution
of the set of equations~(\ref{gj1}), (\ref{gj2}).
It is readily seen that~(\ref{gj1}) holds and
$$
\bigodot_{i=K-k+2}^K f_j^{(i)} \left/ \bigodot_{i=K-k+1}^K f_j^{(i)} \right. = F_j^{(K-k+2)}/F_j^{(K-k+1)}=g_j^{(k)}
$$
which proves~(\ref{gj2}).\CQFD


\paragraph{\bf Proof of Proposition \ref{prop:liebFre}} 

The expression~\eqref{eq:LF_epxression} of the copula $\CLF $ is a simple consequence of the definition of the partition $\mathcal{A}_0,\ldots,\mathcal{A}_K$. We only derive the singular component expression in the case of a product of two terms, $K=2$; the general case follows similar lines. Denote $p_1$ and $q_1$ by $p$ and $q$, then $\CLF (u,v)$ can be respectively written on $\mathcal{A}_0$, $\mathcal{A}_1$, and $\mathcal{A}_2$ by $u$, $u^{1-p}v^q$ and $v$. The cross derivative $\frac{\partial^2\CLF }{\partial u\partial v}(x,y)$ then vanishes on $\mathcal{A}_0$ and $\mathcal{A}_2$, and is equal to $(1-p)qx^{-p}y^{q-1}$ on $\mathcal{A}_1$. Using the formula \citep[Eq. (2.4.1)]{Nelsen},
\begin{equation*}
    \ALF (u,v)=\int_0^u\int_0^v \frac{\partial^2\CLF }{\partial u\partial v}(x,y)\ddr x\ddr y,
\end{equation*}
and dividing the double integral above into the three sets of $[0,u]\times[0,v]$ intersected with $\mathcal{A}_0$, $\mathcal{A}_1$, and $\mathcal{A}_2$ yields 
\begin{equation*}
    \ALF (u,v)=\int_{[0,u]\times[0,v] \cap \mathcal{A}_1} \frac{\partial^2\CLF }{\partial u\partial v}(x,y)\ddr x\ddr y,
\end{equation*}
and therefore routine calculations yield the following expressions:
\begin{equation*}
    \ALF (u,v) = \left\{
    \begin{array}{cc}
         -(1-q)u^{\frac{1-p}{1-q}}+(1-p)u\quad &\text{if} \quad (u,v)\in\mathcal{A}_0,  \\
         -(1-q)u^{\frac{1-p}{1-q}}+u^{1-p}v^q-pv^{\frac{p}{q}}\quad &\text{if} \quad (u,v)\in\mathcal{A}_1,  \\
         qv-pv^{\frac{p}{q}}\quad &\text{if} \quad (u,v)\in\mathcal{A}_2.
    \end{array}
    \right.
\end{equation*}
Taking the difference $\SLF =\CLF -\ALF $ yields a simple expression for the singular component
\begin{equation}\label{eq:LF_singular_K=2}
    \SLF (u,v)=p\min(u, v^{\frac{q}{p}})+(1-q)\min(u^{\frac{1-p}{1-q}}, v)=
    p\min(u^p, v^{q})^{\frac{1}{p}}+(1-q)\min(u^{1-p}, v^{1-q})^{\frac{1}{1-q}}.
\end{equation}
Given that $r_1\leq r_2$ yields $p\leq q$ in the case $K=2$, thus leading to expression~\eqref{eq:LF_singular} when $K=2$. The general case $K>2$ is derived similarly. 
\CQFD


\paragraph{\bf Proof of Proposition \ref{prop3coeffs}}  

(i) Taking account of the identity $\min(x, y)+\max(x, y) = x+y$ yields
\begin{align*}
    \CLF \left(\frac{1}{2},\frac{1}{2}\right)=\prod_{k=1}^K\min(2^{-p_k}, 2^{-q_k}) = 2^{-\sum_{k=1}^K \max(p_k, q_k)}= 2^{\sum_{k=1}^K \min(p_k, q_k)-2}
\end{align*}
and consequently
$\beta(\CLF)=2^{\sum_{k=1}^K \min(p_k, q_k)}-1$.\\
(ii) We use the following expression of Kendall's $\tau$, which is convenient since the copula $\CLF$ may have non null singular component (\citep{Nelsen}, Eq. 5.1.12):
$$
\tau(\CLF)=1-4\int_{[0,1]^2} \frac{\partial}{\partial u} \CLF(u,v)\frac{\partial}{\partial v} \CLF(u,v) \ddr u\ddr v 
= 1-4 \sum_{k=0}^K \int_{{\cal A}_k} \frac{\partial}{\partial u} \CLF(u,v)\frac{\partial}{\partial v} \CLF(u,v) \ddr u\ddr v.
$$
Let $k\in\{0,\dots,K\}$. Then, for all $(u,v)\in {\cal A}_k$, one has
$$
\frac{\partial \CLF}{\partial u}(u,v)= (1-\bar p_k) u^{-\bar p_k}v^{\bar q_k}\; \mbox{ and }\;
\frac{\partial \CLF}{\partial v}(u,v)= \bar q_k u^{1-\bar p_k}v^{\bar q_k-1}.
$$
Besides, remarking that $\bar q_0=0$ and $\bar p_K=1$ shows that both terms $k=0$ and $k=K$ do not contribute
to the sum in $\tau(\CLF)$.
The result 
$$
    \tau(\CLF ) = 1 - \sum_{k=1}^{K-1} \frac{(1-\bar p_k) \bar q_k (r_{k+1}-r_k)}{(\bar q_k r_k + (1-\bar p_k))(\bar q_k r_{k+1} + (1-\bar p_k))}
$$
then follows. \\
(iii) Recall that 
\begin{align*}
\rho(\CLF)
&= 12 \sum_{k=0}^K \int_{{\cal A}_k}\CLF(u,v)\,\ddr u\ddr v-3
=12 \sum_{k=0}^K\int_{{\cal A}_k} u^{1-\bar p_k}  v^{\bar q_k}\,\ddr u\ddr v-3 \\
&= \frac{12(1+r_1+r_1 r_K)}{(2+r_1)(1+2r_K)} -3 + 
    \sum_{k=1}^{K-1} \frac{r_{k+1}-r_k}{((1+\bar q_k) r_k + (2-\bar p_k))((1+\bar q_k) r_{k+1} + (2-\bar p_k))},
  \end{align*}
and the result is then established.
\CQFD

\section{Auxiliary results\label{appendix:aux}}

\begin{Lem}
 \label{lemhomo}
 Let $\varphi:{\mathbb R}_-^d\to {\mathbb R}_-$ be a continuous function such that
 \begin{equation}
  \label{eqhomo}
  \varphi(x)=\varphi(ax) + \varphi((1-a)x)
 \end{equation}
  for all $a\in(0,1)$ and $x\in {\mathbb R}_-^d$.
 Then, necessarily, $\varphi$ is homogeneous of degree 1.
\end{Lem}

\paragraph{\bf Proof} 
Firstly, let us prove by induction the property $(P_n)$: $\varphi(x/n)=\varphi(x)/n$ for all
$n\in{\mathbb N}\setminus\{0\}$ and $x\in {\mathbb R}_-^d$. $(P_1)$ is straightfowardly true. Assume $(P_n)$ holds.
Then,
$$
\varphi\left(\frac{x}{n+1}\right)=\varphi\left(\frac{nx}{n+1} \frac{1}{n}\right) = 
\frac{1}{n}\varphi\left(\frac{nx}{n+1} \right)
$$
and (\ref{eqhomo}) entails
$$
\varphi(x) = \varphi\left(\frac{x}{n+1} \right) + \varphi\left(\frac{nx}{n+1} \right)
= \varphi\left(\frac{x}{n+1} \right) + n \varphi\left(\frac{x}{n+1} \right) = (n+1)  \varphi\left(\frac{x}{n+1} \right),
$$
which proves $(P_{n+1})$.\\
Second, for all $m\in{\mathbb N}\setminus\{0\}$, $(P_m)$ shows that $\varphi(x)=m\varphi(x/m)$ and thus, letting
$y=x/m$, $\varphi(my)=m\varphi(y)$ for all $y\in {\mathbb R}_-^d$. This property can be extended to $m=0$ since
letting $a\to0$ in~(\ref{eqhomo}) yields $\varphi(0)=0$.\\
Third, let $q\in {\mathbb Q}_+$. There exists $(m,n)\in {\mathbb N}\times {\mathbb N}\setminus\{0\}$ such that
$q=m/n$. From the first two points, $\varphi(qx)=\varphi(mx/n)=m\varphi(x/n)=m\varphi(x)/n=q\varphi(x)$.\\
Finally, the continuity of $\varphi$ and the density of ${\mathbb Q}_+$ in ${\mathbb R}_+$ imply that
$\varphi(tx)=t\varphi(x)$ for all $t\in {\mathbb R}_+$ and $x\in {\mathbb R}_-^d$. The result is thus proved.\CQFD


\begin{Lem}\label{theocop}
For all $k\geq 1$ let $C_k$ be a $d$-variate copula and $f^{(k)}_j \in {\cal F}$ for all $j\in\{1,\dots,d\}$, with the assumption $f^{(1)}_j(t)=1$ for all $t\in[0,1]$. 
The sequence  $(\tilde{C}^{(k)})_{k\geq 1}$ defined iteratively by~(\ref{def1}) and~(\ref{def2}) 
is a sequence of $d$-variate copulas.
\end{Lem}

\paragraph{\bf Proof}  The proof is done by induction on $k$. Let ${\cal C}$ be the set of all $d$-variate copulas.
First, it is clear that $\tilde{C}^{(1)}\in {\cal C}$ from~(\ref{def1}). Second, let us assume that 
$\tilde{C}^{(k-1)}\in {\cal C}$ and prove that $\tilde{C}^{(k)}\in {\cal C}$. for all $k\geq 2$.
 Let $(Y_1,\dots,Y_d)$ and $(Z_1,\dots,Z_d)$ be two independent random vectors in $[0,1]^d$ drawn 
respectively from the cdf $\tilde{C}^{(k-1)}(f_1^{(k)},\dots,f_d^{(k)})$ and 
${C}_K(\Id/f_1^{(k)},\dots,\Id/f_d^{(k)})$. For all $j\in\{1,\dots,d\}$ define the random variable 
$X_j=\max(Y_j,Z_j)$. For all $\bm{u}\in[0,1]^d$, the cdf of $(X_1,\dots,X_d)$ is given by
\begin{align*}
{\Pr}(X_1\leq u_1,\dots,X_d\leq u_d) &= {\Pr}(Y_1\leq u_1,\dots,Y_d\leq u_d)
{\Pr}(Z_1\leq u_1,\dots,Z_d\leq u_d)\\
&= \tilde{C}^{(k-1)}(f_1^{(k)}(u_1),\dots,f_d^{(k)}(u_d)) {C}_K(u_1/f_1^{(k)}(u_1),\dots,x_d/f_d^{(k)}(u_d)) 
=\tilde{C}^{(k)}(\bm{u}),
\end{align*}
from (\ref{def2}). This proves that $\tilde{C}^{(k)}$ is a cdf. The margins are uniform by construction. \CQFD


\bibliographystyle{model2-names}
\bibliography{biblio}

\end{document}